\documentclass[times,final]{article}

\usepackage{fullpage}
\usepackage{framed,multirow}

\usepackage{latexsym}

\usepackage{amsmath, amssymb, mathabx, mathtools}
\usepackage{graphicx, psfrag}
\usepackage{empheq}
\usepackage{subcaption}
\usepackage{bm}
\usepackage{xcolor}
\usepackage[capitalize]{cleveref}
\usepackage[numbers]{natbib}
\usepackage{authblk}

\graphicspath{{./images/}}

\def\mesh{\mathcal{T}}
\def\dualmesh{{\widetilde{\mathcal{T}}}}
\def\trig{T}
\def\dualtrig{\widetilde {T}}

\def\curl{\operatorname{curl}}                                
\usepackage{colortbl}

\usepackage{url}

\title{A higher-order dual cell method for time-domain Maxwell equations}
\author[1]{Lorenzo Codecasa}
\author[2]{Bernard Kapidani}
\author[3]{Joachim Sch\"oberl}
\author[3]{Markus Wess}

\affil[3]{Institute of Analysis and Scientific Computing, Technische Universit\unexpanded{\"a}t Wien, A-1040, Vienna, Austria.\authorcr
    \tt markus.wess@tuwien.ac.at, joachim.schoeberl@tuwien.ac.at }
\affil[2]{Institute of Mathematics,  \unexpanded{\'{E}}cole Polytechnique F\unexpanded{\'{e}}d\unexpanded{\'{e}}rale Lausanne, CH-1015 Lausanne, Switzerland \authorcr
        \tt bernard.kapidani@epfl.ch }
\affil[1]{Dipartimento di Elettronica, Informatica e Bioingegneria, Politecnico di Milano, I-20133 Milano, Italy \authorcr
    \tt lorenzo.codecasa@polimi.it }

\date{\today}

\begin{document}
\maketitle
\begin{abstract}
{
We present a higher-order extension of the dual cell method for the time-domain Maxwell equations in three spatial dimensions. The approach builds upon a variational reinterpretation of the Finite Integration Technique on dual meshes and generalises a previously developed two-dimensional high-order formulation. The electric and magnetic fields are discretised on mutually dual barycentric grids using curl-conforming polynomial spaces constructed via tensor-product Gauss--Radau interpolation. The resulting semi-discrete formulation yields block-diagonal mass matrices and sparse discrete curl operators, enabling explicit time integration while preserving a discrete energy identity. 
Special attention is devoted to the construction of compatible approximation spaces on the three-dimensional primal and dual meshes, the reference-to-physical element mappings, and the preservation of tangential continuity. We show that the method achieves arbitrary-order convergence, avoids spurious modes, and maintains optimal sparsity properties. Numerical experiments confirm spectral correctness, high-order accuracy, and computational efficiency on unstructured tetrahedral meshes.
}
\end{abstract}

\section{Introduction}\label{sec:intro}

The numerical approximation of time-dependent initial boundary value problems for hyperbolic systems such as Maxwell's equations is traditionally dominated by finite difference schemes, most prominently the second-order staggered grid formulation introduced by Yee~\cite{yeeNumericalSolutionInitial1966}, commonly referred to as the Finite Difference Time Domain (FDTD) method. The widespread adoption of FDTD is largely due to its algorithmic simplicity, structured stencil, and suitability for large-scale parallel implementations~\cite{tafloveAdvancesFDTDComputational2013}.
In many engineering applications, however, electromagnetic devices exhibit complex geometrical features that are more naturally described by unstructured tetrahedral meshes. This motivates the use of finite element methods (FEM) for the time-domain Maxwell equations~\cite{jinFiniteElementMethod2014a}. In particular, edge elements introduced by N\'ed\'elec~\cite{nedelecMixedFiniteElements1980} and their high-order extensions~\cite{schoberlHighOrderNedelec2005} provide curl-conforming discretisations that eliminate spurious gradient modes at the continuous level.

A major limitation of standard continuous FEM discretisations in explicit time-domain simulations lies in the presence of non-diagonal mass matrices. The mass matrix associated with edge element spaces is sparse but not diagonal, and its inversion is computationally demanding at every time step. While iterative solvers such as preconditioned conjugate gradient methods may be employed, this increases computational cost and may slightly deteriorate the exact energy conservation property at the semi-discrete level.

Discontinuous Galerkin (DG) methods provide an alternative that restores locality and parallel scalability by allowing element-wise discontinuities~\cite{hesthavenNodalDiscontinuousGalerkin2008,arnoldUnifiedAnalysisDiscontinuous2002}. However, DG schemes require carefully designed numerical fluxes and stabilisation mechanisms, and may exhibit spurious solutions or instabilities unless stabilisation parameters are appropriately tuned.
The present work extends a line of research initiated by the authors, rooted in the Finite Integration Technique (FIT)~\cite{weilandTimeDomainElectromagnetic1996,codecasaExplicitConsistentConditionally2008,codecasaNovelFDTDTechnique2018}. In~\cite{kapidaniTimeDomainCellMethod2020}, it was shown that the low-order FIT formulation can be interpreted as a variational method on barycentric dual meshes. 
A numerical analysis of closely related lowest order methods is presented in  ~\cite{BrezziFortinMarini, LeeWinther,EggerRaduYee}.

In~\cite{kapidaniArbitraryorderCellMethod2021,WessKapidaniCodecasaSchoeberl}, a high-order generalisation was developed in two dimensions using orthogonal polynomial bases to achieve arbitrary-order convergence together with optimally sparse mass matrices.
Existing methods based on Discrete Exterior Calculus ~\cite{hiranianilnirmalDiscreteExteriorCalculus2003,hiraniDelaunayHodgeStar2013} provide lowest order explicit schemes for hyperbolic systems of equations on unstructured dual grids and more recent research has been published to provide mass lumping schemes for scalar wave equations and Maxwell equations in three dimensions up to second order space accuracy~\cite{geeversNewHigherOrderMassLumped2018,eggerMasslumpedMixedFinite2020,eggerSecondOrderFiniteElement2021}, but no general recipe is available to extend these methods to arbitrary order of convergence to render them competitive with full DG schemes in terms of performance.
In this work, we extend the high-order dual cell method to three spatial dimensions and unstructured grids (similar ideas for staggered high-order schemes on cuboidal meshes have been investigated, e.g., in~\cite{chungConvergenceSuperconvergenceStaggered2013,chungStaggeredDGMethod2014}). As in the two-dimensional case, the electric and magnetic fields are discretised on mutually dual barycentric grids. One field is defined on the primal tetrahedral mesh, while the other is discretised on the associated polyhedral dual mesh obtained by barycentric subdivision. A central objective of the present formulation is to retain the structural properties of Maxwell's equations at the discrete level, in particular the absence of spurious modes and the conservation of electromagnetic energy, while at the same time achieving favourable sparsity and explicit time stepping through block-diagonal (lumped) mass matrices.
The resulting scheme preserves the key properties of the two-dimensional formulation: arbitrary high-order convergence, freedom from spurious solutions, explicit time stepping, and a discrete energy conservation law. The extension to three dimensions significantly broadens the applicability of the method to realistic electromagnetic configurations.

{The paper is structured as follows. After recalling the main ingredients of the weak formulation and the construction of barycentric dual meshes in \cref{sec:formulation}, we introduce the three-dimensional extension in \cref{sec:spacedisc}, detailing the approximation spaces, local mappings, and discrete operators. We then describe our choice of time discretisation in \cref{sec:timedisc}. The numerical behaviour of the method is investigated in \cref{sec:numerics}, where convergence studies and benchmark tests confirm arbitrary-order convergence, spectral correctness, and computational efficiency of the extension to three dimensions. We conclude with a short summary and an outlook on future research directions.}

\section{Weak Formulation {on Dual Grids}}\label{sec:formulation}

{Let $\Omega \subset \mathbb{R}^3$ be a bounded Lipschitz domain.}
The time-domain Maxwell equations in $\Omega$ can be written in the form
\begin{align*}
& \boldsymbol{\varepsilon}\,\partial_t \bm{E} = \curl \bm{H},\\
& \boldsymbol{\mu}\,\partial_t \bm{H} = -\curl \bm{E},
\end{align*}
where the electric field $\bm{E}$ and the magnetic field $\bm{H}$ are functions of the position vector $\mathbf{x}=(x_1, x_2, x_3)$ and time instant $t$. 
The electric permittivity $\boldsymbol{\varepsilon}$  and the magnetic permeability $\boldsymbol{\mu}$  are symmetric positive definite double tensors {(i.e., symmetric positive definite matrix-valued functions)} and may depend on the position $\mathbf{x}$.
For the sake of simplicity the boundary $\partial \Omega$ is assumed to be an electric wall
{(i.e., $\bm{n}\times\bm{E}=0$ on $\partial\Omega$)}.

{We assume throughout that $\boldsymbol{\varepsilon},\boldsymbol{\mu}\in L^\infty(\Omega)^{3\times3}$ and are uniformly positive definite.}

\begin{figure}
  \centering
  \begin{subfigure}{0.24\textwidth}
    \centering
    \includegraphics[scale=0.5]{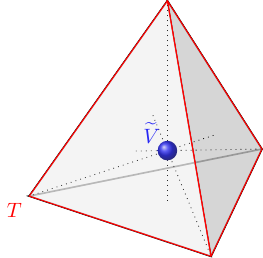}
    \caption{Primal cell and dual vertex}
    \label{fig:barycentric_v}
  \end{subfigure}
  \begin{subfigure}{0.24\textwidth}
    \centering
    \includegraphics[scale = 0.5]{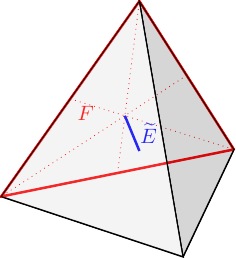}
    \caption{Primal face and dual edge}
    \label{fig:barycentric_e}
  \end{subfigure}
  \begin{subfigure}{0.24\textwidth}
    \centering
    \includegraphics[scale = 0.5]{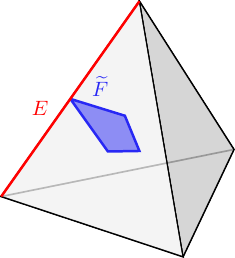}
    \caption{Primal edge and dual face}
    \label{fig:barycentric_f}
  \end{subfigure}
  \begin{subfigure}{0.24\textwidth}
    \centering
    \includegraphics[scale = 0.5]{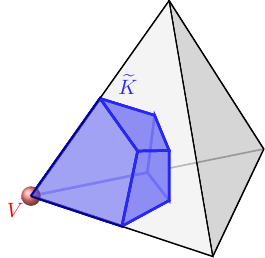}
    \caption{Primal vertex and dual cell}
    \label{fig:barycentric_c}
  \end{subfigure}
  \caption{Entities of the primal mesh (red) and their dual counterparts (blue) for a mesh consisting of one single primal cell.}
  \label{fig:barycentric}
\end{figure}

Let $\mesh$ be a primal simplicial complex covering $\Omega$, 
 and let $\dualmesh$ be its dual, obtained by barycentric subdivision of  $\mesh$ as follows: 
\begin{itemize}
  \item Each vertex $\widetilde{V}$ of $\dualmesh$ is the dual of one cell $\trig$ of $\mesh$,  located in the barycentre of $\trig$ (see \cref{fig:barycentric_v}).

  \item Each edge $\widetilde{E}$ of $\dualmesh$ is the dual of one face $F$ of $\mesh$. In each cell $\trig$ of $\mesh$, dual of vertex $\widetilde{V}$ of  $\dualmesh$, which has face $F$ in its boundary,
$\widetilde{E}$ is the segment joining $\widetilde{V}$ with the barycentre of $F$ (see \cref{fig:barycentric_e}).

\item Each face $\widetilde{F}$ of $\dualmesh$ is the dual of one edge $E$ of $\mesh$. 
Let $\trig$ be any cell of $\mesh$, dual of vertex $\widetilde{V}$ of  $\dualmesh$, having $E$ on its boundary, and let 
$F_1, F_2$ be the two faces of $\mesh$, dual of edges $\widetilde{E}_1$, $\widetilde{E}_2$ of $\dualmesh$ respectively,  belonging to $\trig$ and containing $E$. 
The restriction of $\widetilde{F}$ to $\trig$ is the quadrilateral having  $\widetilde{E}_1$, $\widetilde{E}_2$ as edges and the barycentre of $E$ as a vertex (see \cref{fig:barycentric_f}).

\item Each cell $\dualtrig$ of $\dualmesh$ is the dual of one node $N$ of $\mesh$. 
Let $\trig$ be any cell of $\mesh$, dual of vertex $\widetilde{V}$ of  $\dualmesh$, having $N$ on its boundary, and let 
$E_1, E_2, E_3$ be the three edges of $\mesh$, dual of faces $\widetilde{F}_1$, $\widetilde{F}_2$, $\widetilde{F}_3$ of $\dualmesh$ respectively,  belonging to $\trig$ and containing $V$. 
The restriction of $\dualtrig$ to $\trig$ is the hexahedron having  $\widetilde{F}_1$, $\widetilde{F}_2$, $\widetilde{F}_3$ as faces and $N$ as a vertex (see \cref{fig:barycentric_c}).
\end{itemize}
{This construction yields a pair of mutually dual cell complexes satisfying the usual incidence relations of discrete differential geometry.}

\begin{figure}
  \centering
  \includegraphics[scale = 0.8]{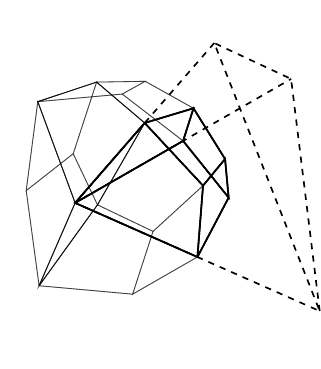}
  \caption{A primal (tetrahedral) and dual (polyhedral) cell.}
\end{figure}

Let $\mathrm{H}(\curl,\mesh)$ be the space of square integrable vector functions whose restriction to each cell ${T}$ of $\mesh$ has square integrable curl.
Similarly let $\mathrm{H}(\curl,\dualmesh)$ be the space of square integrable vector functions whose restriction to each cell ${\widetilde{T}}$ of $\dualmesh$ has square integrable curl.
{More precisely, these are broken Sobolev spaces, i.e.,
\[
\mathrm{H}(\curl,\mesh)
=
\{ \bm{v}\in L^2(\Omega)^3 \;:\; \bm{v}|_T \in H(\curl,T)\ \forall T\in\mesh \},
\]
and analogously for $\dualmesh$.}
By applying Stokes' theorem, the following weak formulation of Maxwell equations ensues. Let $\bm{E}$  and $\bm{H}$ be square integrable vector functions with square integrable curl.  For each  $\bm{e}\in \mathrm{H}(\curl,\dualmesh)$ it results in
\begin{align}\label{E1}
& \sum_{\widetilde{T} \in \dualmesh} \int_{\widetilde{T}} \bm{e} \cdot \boldsymbol{\varepsilon}\partial_t \bm{E} \,d V =
 \sum_{\widetilde{T} \in \dualmesh} \left(-\int_{\partial \widetilde{T}\backslash\partial \Omega} \bm{e} \times \bm{H} \cdot \mathbf{n}\, d S + 
 \int_{\widetilde{T}} \curl\bm{e} \cdot \bm{H}\, d V\right),
\end{align} 
 where $\mathbf{n}$ is the unit vector function outward normal to $\partial \widetilde{T}$.

Similarly for each
 $\bm{h}\in \mathrm{H}(\curl,\mesh)$ we obtain
  \begin{align}\label{E2}
 & \sum_{T \in \mesh} \int_{T} \bm{h} \cdot \boldsymbol{\mu}\partial_t \bm{H} \,d V = 
- \sum_{T \in \mesh} \left(-\int_{\partial T} \bm{h} \times \bm{E} \cdot \mathbf{n}\, d S +
 \int_{T} \curl\bm{h} \cdot \bm{E}\, d V\right).
\end{align} 

{Let now $\mathcal{K}$ be the set of all non-empty intersections $K = T \cap \widetilde{T}$, where $T$ is a cell of the primal mesh $\mathcal{T}$ and $\widetilde{T}$ is a cell of the dual mesh $\widetilde{\mathcal{T}}$. For each $K \in \mathcal{K}$ we further define the boundary portions $\partial_{\mathcal{T}} K = \partial K \cap \partial T$ and $\partial_{\widetilde{\mathcal{T}}} K = \partial K \cap \partial \widetilde{T}$.}

\begin{figure}
  \centering
  \includegraphics[scale = 0.8]{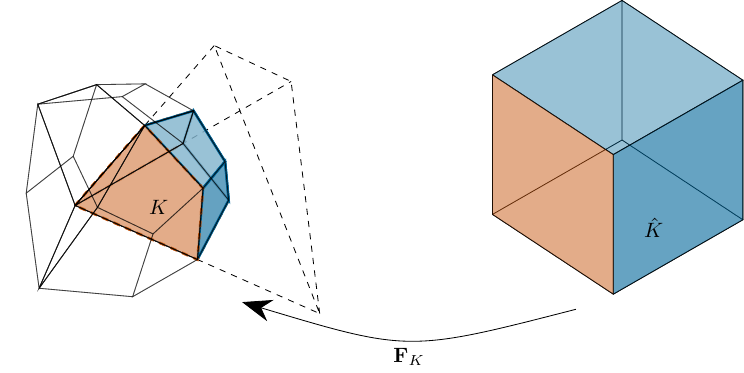}
  \caption{{Reference-to-physical trilinear mapping $\mathbf{F}_K$ from the cube $\widehat{K}$ to a generic polyhedral subcell $K\in\mathcal{K}$.}}
\end{figure}

Eqs. (\ref{E1}), (\ref{E2}) can then be rewritten in the form
 \begin{align}
&  \sum_{K \in \mathcal{K}} \int_{K} \bm{e} \cdot \boldsymbol{\varepsilon}\partial_t \bm{E}\, d V= 
 \sum_{K \in \mathcal{K}} \left(-\int_{\partial K_{\dualmesh}} \bm{e} \times \bm{H} \cdot \mathbf{n}\, d S+ 
 \int_{K} \curl\bm{e} \cdot \bm{H}\, d V\right), \label{E3}
 \end{align} 
 \begin{align}
  \sum_{K \in \mathcal{K}} \int_{K} \bm{h} \cdot \boldsymbol{\mu}\partial_t \bm{H} \,d V = 
- \sum_{K \in \mathcal{K}} \left(-\int_{\partial K_{\mesh}} \bm{h} \times \bm{E} \cdot \mathbf{n}\, d S +
 \int_{K} \curl\bm{h} \cdot \bm{E}\, d V\right). \label{E4}
\end{align}

Let the position vector  $\mathbf{{x}}$ in each hexahedron $K$ of $\mathcal{K}$ be a trilinear map $\mathbf{x}=\mathbf{F}_{K}(\widehat{\mathbf{x}})$ of the position vector $\widehat{\mathbf{x}}=(\widehat{x}_1, \widehat{x}_2, \widehat{x}_3)$ in the cube $\widehat{K}=[-1, 1]^3$. 
Let the vertex of $K$ corresponding to the node of $\mesh$ be the map of  $\widehat{\mathbf{x}}=(-1,-1,-1)$, and the vertex of $K$ corresponding to the node of $\dualmesh$ be the map of  $\widehat{\mathbf{x}}=(1,1,1)$.
 Let $\mathbf{J}_{K}(\widehat{\mathbf{x}})=\partial \mathbf{F}_K(\widehat{\mathbf{x}})/\partial\widehat{\mathbf{x}}$ be the Jacobean matrix of this map.
{Vector fields on $K$ are related to reference fields on $\widehat{K}$ via the covariant Piola transformation
\[
\bm{v} = \mathbf{J}_K^{-T} \widehat{\bm{v}},
\]
which preserves tangential continuity across element interfaces and yields curl-conforming approximation spaces.}
Eqs. (\ref{E3}), (\ref{E4}) can then be rewritten in the form
\begin{align}
  \sum_{K \in \mathcal{K}} \int_{\widehat{K}} \widehat{\bm{e}}_{\widehat{K}} \cdot \widehat{\boldsymbol{\varepsilon}}_K\partial_t \widehat{\bm{E}}_{\widehat{K}} \, |\mathbf{J}_{K}|d\widehat{ V} =
\sum_{K \in \mathcal{K}} \left(-\int_{\partial \widehat{K}_{\dualmesh}} \widehat{\bm{e}}_{\widehat{K}} \times \widehat{\bm{H}}_{\widehat{K}} \cdot \widehat{\mathbf{n}}\, d\widehat{ S} + 
  \int_{\widehat{K}}{\widehat{\curl}\, \widehat{\bm{e}}_{\widehat{K}}}\cdot \widehat{\bm{H}}_{\widehat{K}}\, d\widehat{ V}\right), \label{E5}
\end{align}
\begin{align}
  \sum_{K \in \mathcal{K}} \int_{\widehat{K}} \widehat{\bm{h}}_{\widehat{K}} \cdot  \widehat{\boldsymbol{\mu}}_K \partial_t \widehat{\bm{H}}_{\widehat{K}} \, |\mathbf{J}_{K}|d\widehat{ V} =
- \sum_{K \in \mathcal{K}} \left(-\int_{\partial \widehat{K}_{\mesh}} \widehat{\bm{h}}_{\widehat{K}} \times \widehat{\bm{E}}_K\cdot \widehat{\mathbf{n}}\, d \widehat{ S} +
  \int_{\widehat{K}} {\widehat{\curl}\, \widehat{\bm{h}}_{\widehat{K}}} \cdot \widehat{\bm{E}}_{\widehat{K}}\, d\widehat{ V}\right), \label{E6}
\end{align}
in which: $\bm{E}= \mathbf{J}_{K}^{-T}\widehat{\bm{E}}_{\widehat{K}}$, $\bm{e}= \mathbf{J}_{K}^{-T}\widehat{\bm{e}}_{\widehat{K}}$, $\bm{H}= \mathbf{J}_{K}^{-T}\widehat{\bm{H}}_{\widehat{K}}$, $\bm{h}= \mathbf{J}_{K}^{-T}\widehat{\bm{h}}_{\widehat{K}}$, $\widehat{\boldsymbol{\varepsilon}}_K=\mathbf{J}_{K}^{-1}\boldsymbol{\varepsilon}\mathbf{J}_{K}^{-T}$ and $\widehat{\boldsymbol{\mu}}_K=\mathbf{J}_{K}^{-1}\boldsymbol{\mu}\mathbf{J}_{K}^{-T}$; $d\widehat{V}$ is an infinitesimal volume, $d\widehat{S}$ is an infinitesimal area; $\widehat{\curl}$ is the curl in $\widehat{K}$; $\widehat{\mathbf{n}}$ is the unit vector function outward normal to $\partial \widehat{K}$; $\partial\widehat{K}_{\mesh}$ and $\partial\widehat{K}_{\dualmesh}$ are mapped by $\mathbf{F}_K$ into $\partial{K}_{\mesh}$ and $\partial{K}_{\dualmesh}$  respectively.

\section{Discretisation}\label{sec:disc}
{In the following we present how the weak formulation \eqref{E5}, \eqref{E6} are discretised in space and time. For the spatial discretisation we construct in subsection \ref{sec:spacedisc} finite dimensional subspaces of $\mathrm{H}(\curl,K)$ for each hexahedron $K$ of arbitrary polynomial order. These subspaces are subsequentially combined to $\mathrm H(\curl)$ conforming sub spaces on primal and dual cells respectively. In subsection \ref{sec:timedisc} we briefly discuss the time discretisation}
\subsection{Spatial Discretisation}\label{sec:spacedisc}
  
Let $\mathcal Q_P(\mathcal{K})$ with $P>0$ be the set of vector functions whose restriction to each hexahedron $K\in\mathcal{K}$ corresponding to cube $\widehat{K}$ is  $\mathbf{v} =\mathbf{J}_{K}^{-T}\widehat{\mathbf{v}}_{{K}}$  being $\widehat{\mathbf{v}}_{{K}}$ a polynomial vector function of degree at most $P$ in each of the variables $\widehat{x}_1$, $\widehat{x}_2$, $\widehat{x}_3$. 

{
Using the trilinear mapping $\mathbf{F}_K : \widehat{K} \to K$ and associated Jacobian $\mathbf{J}_K$ introduced in Section~\ref{sec:formulation}, we define the tensor-product polynomial space
\[
\mathcal{Q}_P(\mathcal{K})
=
\left\{
\bm{v} \; : \;
\bm{v}|_K = \mathbf{J}_K^{-T} \widehat{\bm{v}}, 
\quad 
\widehat{\bm{v}} \in \mathbb{Q}_P^3(\widehat{K})
\right\},
\]
where $\mathbb{Q}_P$ denotes polynomials of degree at most $P$ in each coordinate.
For each $K$ we employ tensor-product Gauss--Radau interpolation polynomials associated with nodes $-1 = \widehat{\xi}_0 < \widehat{\xi}_1 < \dots < \widehat{\xi}_P < 1$ on $[-1,1]$ in every coordinate direction, which will be used both for defining basis functions and for constructing high-order accurate quadrature rules.}

{Curved geometries can be treated in the usual isoparametric fashion by choosing $\mathbf{F}_K$ as a higher-order polynomial (or spline) mapping that approximates curved boundaries. Provided that $\mathbf{F}_K$ is sufficiently smooth and of polynomial degree at least $P$, the induced Piola mapping preserves the convergence order of the approximation in physical space, in complete analogy with standard high-order isoparametric finite elements~\cite{jinFiniteElementMethod2014a}.
}

\subsubsection{$H_P(\curl, \widetilde{\mesh})$}
Let it be $H_P(\curl, \widetilde{\mesh})=H(\text{curl},\widetilde{ \mesh}) \bigcap \mathcal Q_P(\mathcal{K})$. 
A basis 
for $H_P(\curl, \dualmesh)$ is the set of  vector functions ${\bm{E}}_n$, with $n=1,\dots, N$, whose 
restriction  to any hexahedron $K$ of $\mathcal{K}$ corresponding to cube $\widehat{K}$ is either identically zero or given by
\begin{align}
&{\bm{E}_{n,K}} =\mathbf{J}_{K}^{-T}\widehat{{\bm{E}}}_{{n,\widehat{K}}},\nonumber\\
&\widehat{{\bm{E}}}_{{n,\widehat{K}}}= \widehat{\ell}_{\boldsymbol{\alpha}}(\widehat{\mathbf{x}})\widehat{\mathbf{e}}_{i}, \label{E7} 
\end{align}
in which $i\in\{1,2,3\}$, $\widehat{\mathbf{e}}_{i}$ is the unit vector corresponding to axis $\widehat{x}_{i}$, $\boldsymbol{\alpha}=(\alpha_1, \alpha_2, \alpha_3)\in\{0,\dots,P\}^3$, $\widehat{\ell}_{\boldsymbol{\alpha}}(\widehat{\mathbf{x}})=\widehat{\ell}_{\alpha_1}(\widehat{x}_1)\cdot  \widehat{\ell}_{\alpha_2}(\widehat{x}_2) \cdot \widehat{\ell}_{\alpha_3}(\widehat{x}_3)$, in which
$\widehat{\ell}_{\alpha}(\widehat{x})$, with $0\leq\alpha\leq P$, is the interpolation polynomial in $\widehat{x}$, equal to one in the interpolation point  $\widehat{\xi}_{\alpha}$, having as interpolation points the $P+1 $ Gauss-Radau quadrature points $-1=\widehat{\xi}_0<\widehat{\xi}_1<\cdots<\widehat{\xi}_P<1$ for interval $[-1, 1]$, corresponding to weights $\widehat{w}_0, \widehat{w}_1, \cdots, \widehat{w}_P$.
Three types of basis vector functions follow:
\\
\\
1. {\bf Cell-centred basis functions.} 
\\
    A permutation $(i\, j\, k)$  of $(1\, 2\, 3)$ exists such that $\alpha_j>0$ and $\alpha_k>0$. In this case ${\bm{E}}_{{n}}$ is zero on the two faces of ${\mathcal{F}}$ belonging to the boundary of $K$ and corresponding to the faces of $\widehat{K}$ having equations $\widehat{x}_j=-1$ and $\widehat{x}_k=-1$. Besides ${\bm{E}}_{{n}}$ is normal to the third face of  ${\mathcal{F}}$ belonging to the boundary of $K$ and corresponding to the face of $\widehat{K}$ having equation $\widehat{x}_i=-1$. As a result ${\bm{E}}_{{n}}$ has zero tangential trace at all faces of ${\mathcal{F}}$ on $\partial K$ and can be uniquely extended {$\mathrm{H}(\curl, \widetilde{\mesh})$ conforming by $\mathbf{0}$ to the all hexahedra $K'\in\mathcal{K}$ with $K'\neq K$. The number of  vector functions of this kind is $3(P+1)P^2|\mathcal{K}|$.
\\
\\
2. {\bf Face-centred basis functions.} 
\\
 A permutation $(i\, j\, k)$  of $(1\, 2\, 3)$ exists such that $\alpha_{j}=0$ and $\alpha_{k}>0$. In this case ${\bm{E}}_{n}$ is zero on the face of ${\mathcal{F}}$ belonging to the boundary of $K$ and corresponding to the face of $\widehat{K}$ having equation $\widehat{x}_k=-1$.
Besides ${\bm{E}}_{{n}}$ is normal to the face of ${\mathcal{F}}$ belonging to the boundary of $K$ and corresponding to the face of $\widehat{K}$ of equation $\widehat{x}_i=-1$. 
In order to be {$\mathrm{H}(\curl,\widetilde{\mesh})$ conforming, ${\bm{E}}_{{n}}$ is uniquely extended to all hexahedra $K'\in\mathcal{K}$ as follows:
if $K'$ contains the face ${F}\in{\mathcal{F}}$ belonging to the boundary of $K$ and corresponding to the face of $\widehat{K}$ of equation $\widehat{x}_{j}=-1$, then  ${\bm{E}}_{{n}}$ is extended to $K'$ by $\widehat{\bm{E}}_{{{n}},\widehat{K}'}=\widehat{\mathbf{e}}_{i'}\widehat{\ell}_{\boldsymbol{\alpha}'}(\widehat{\mathbf{x}})$, $\boldsymbol{\alpha}'$ being a permutation of $\boldsymbol{\alpha}$;
otherwise ${\bm{E}}_{{n}}$ is extended by $\mathbf{0}$ to $K'$.
  The number of  vector functions of this kind is $2(P+1)P|{\mathcal{F}}|$.
\\
\\
3. {\bf Edge-centred basis functions. }
\\
A permutation $(i\, j\, k)$  of $(1\, 2\, 3)$ exists such that  $\alpha_{j}=0$ and $\alpha_{k}=0$.  In this case ${\bm{E}}_{{n}}$ is normal to the face of ${\mathcal{F}}$ belonging to the boundary of $K$ and corresponding to the face of $\widehat{K}$ of equation $\widehat{x}_i=-1$. 
In order to be {$\mathrm{H}(\curl, \widetilde{\mesh})$ conforming, ${\bm{E}}_{{n}}$ is uniquely extended to all hexahedra $K'\in\mathcal{K}$ as follows:
if $K'$ contains the edge ${E}\in{\mathcal{E}}$ belonging to the boundary of $K$ and corresponding to the edge of $\widehat{K}$ of equations $\widehat{x}_{j}=-1$, $\widehat{x}_{k}=-1$, then  ${\bm{E}}_{{n}}$ is extended to $K'$ by $\widehat{\bm{E}}_{{{n}},\widehat{K}'}=\widehat{\mathbf{e}}_{i'}\widehat{\ell}_{\boldsymbol{\alpha}'}(\widehat{\mathbf{x}})$, $\boldsymbol{\alpha}'$ being a permutation of $\boldsymbol{\alpha}$;
otherwise ${\bm{E}}_{{n}}$ is extended by $\mathbf{0}$ to $K'$.
The number of  vector functions of this kind is $(P+1) |{\mathcal{E}}|$.
\\

It is noticed that the support of each vector function $\bm{E}_{{n}}$ in all these three cases is contained in a single cell $\widetilde{T}\in\widetilde{\mesh}$ and consists of one, two or three hexahedra. 

\begin{figure}
  \centering
  \begin{subfigure}{0.3\textwidth}
    \centering
    \includegraphics[scale=0.4]{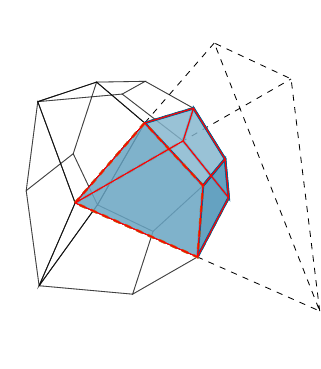}
    \caption{Support of a cell-centered basis function}
  \end{subfigure}
  \begin{subfigure}{0.3\textwidth}
    \centering
    \includegraphics[scale=0.4]{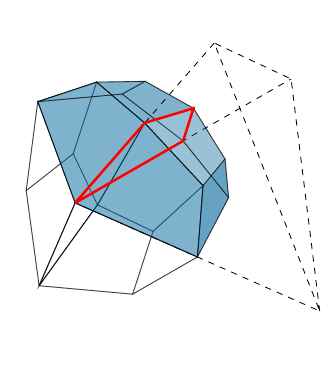}
    \caption{Support of a face-centered basis function}
  \end{subfigure}
  \begin{subfigure}{0.3\textwidth}
    \centering
    \includegraphics[scale=0.4]{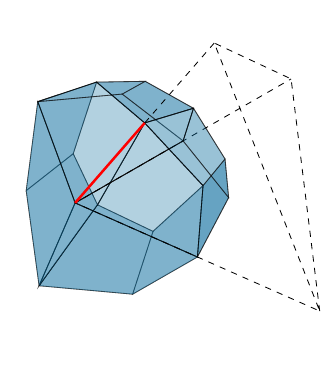}
    \caption{Support of a edge-centered basis function}
  \end{subfigure}
  \caption{Support of different types of basis functions on the dual mesh}
\end{figure}
\begin{figure}
  \centering
  \begin{subfigure}{0.3\textwidth}
    \centering
    \includegraphics[scale=0.4]{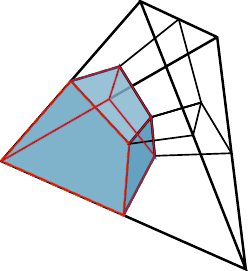}
    \caption{Support of a cell-centered basis function}
  \end{subfigure}
  \begin{subfigure}{0.3\textwidth}
    \centering
    \includegraphics[scale=0.4]{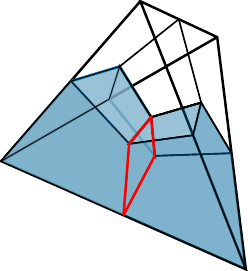}
    \caption{Support of a face-centered basis function}
  \end{subfigure}
  \begin{subfigure}{0.3\textwidth}
    \centering
    \includegraphics[scale=0.4]{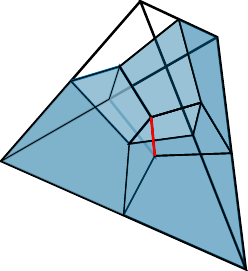}
    \caption{Support of a edge-centered basis function}
  \end{subfigure}
  \caption{Support of different types of basis functions on the primal mesh}
\end{figure}

\subsubsection{$H_P(\curl, {\mesh})$}
In a similar way we define $H_P(\curl, {\mesh})=H(\text{curl},\mathcal{ {T}}) \bigcap \mathcal Q_P(\mathcal{K})$. 
A basis 
 for $H_P(\curl, \mathcal{{T}})$ is the set of  vector functions ${\bm{H}}_{\widetilde{n}}$, with $\widetilde{n}=1,\dots,\widetilde{N}$, whose 
restriction to any hexahedron $K$ of $\mathcal{K}$ corresponding to cube $\widehat{K}$ is either identically zero or given by
\begin{align}
& {\bm{H}}_{\widetilde{n}}=\mathbf{J}_{K}^{-T}\widehat{{\bm{H}}}_{{{\widetilde{n}},\widehat{K}}},\nonumber \\
&\widehat{{\bm{H}}}_{{\widetilde{n}},\widehat{K}}= \widehat{\ell}_{\boldsymbol{\alpha}}(-\widehat{\mathbf{x}}) \widehat{\mathbf{e}}_{i}.\label{E8}
\end{align}
in which $\boldsymbol{\alpha}\in\{0,\dots,P\}^3$.
Three types of basis vector functions follow:
\\
\\
1. {\bf Cell-centred basis functions. }
\\
A permutation $(i\, j\, k)$  of $(1\, 2\, 3)$ exists such that $\alpha_j>0$ and $\alpha_k>0$. In this case ${\bm{H}}_{\widetilde{n}}$ is zero on the two faces of $\widetilde{\mathcal{F}}$ belonging to the boundary of $K$ and corresponding to the faces of $\widehat{K}$ having equations $\widehat{x}_j=1$ and $\widehat{x}_k=1$. Besides ${\bm{H}}_{\widetilde{n}}$ is normal to the third face of  $\widetilde{\mathcal{F}}$ belonging to the boundary of $K$ and corresponding to the face of $\widehat{K}$ having equation $\widehat{x}_i=1$. As a result ${\bm{H}}_{\widetilde{n}}$ has zero tangential trace at all faces of $\widetilde{\mathcal{F}}$ on $\partial K$ and can be uniquely extended {$\mathrm{H}(\curl, {\mesh})$ conforming by $\mathbf{0}$ to  all hexahedra $K'\in\mathcal{K}$ with $K'\neq K$. The number of  vector functions of this kind is $3(P+1)P^2|\mathcal{K}|$.
\\
\\
2. {\bf Face-centred basis functions.}
\\
 A permutation $(i\, j\, k)$  of $(1\, 2\, 3)$ exists such that $\alpha_{j}=0$ and $\alpha_{k}>0$. In this case ${\bm{H}}_{\widetilde{n}}$ is zero on the face of $\widetilde{\mathcal{F}}$ belonging to the boundary of $K$ and corresponding to the face of $\widehat{K}$ having equation $\widehat{x}_k=1$.
Besides ${\bm{H}}_{\widetilde{n}}$ is normal to the face of $\widetilde{\mathcal{F}}$ belonging to the boundary of $K$ and corresponding to the face of $\widehat{K}$ of equation $\widehat{x}_i=1$. 
In order to be {$\mathrm{H}(\curl, {\mesh})$ conforming, ${\bm{H}}_{\widetilde{n}}$ is uniquely extended to all hexahedra $K'\in\mathcal{K}$ as follows:
if $K'$ contains the face $\widetilde{F}\in\widetilde{\mathcal{F}}$ belonging to the boundary of $K$ and corresponding to the face of $\widehat{K}$ of equation $\widehat{x}_{j}=1$, then  ${\bm{H}}_{\widetilde{n}}$ is extended to $K'$ by $\widehat{\bm{H}}_{{\widetilde{n}},\widehat{K}'}=\widehat{\mathbf{e}}_{i'}\widehat{\ell}_{\boldsymbol{\alpha}'}(-\widehat{\mathbf{x}})$, $\boldsymbol{\alpha}'$ being a permutation of $\boldsymbol{\alpha}$;
otherwise ${\bm{H}}_{\widetilde{n}}$ is extended by $\mathbf{0}$ to $K'$.
  The number of  vector functions of this kind is $2(P+1)P|\widetilde{\mathcal{F}}|$.
 \\
\\
3.  {\bf Edge-centred basis functions.}
\\
 A permutation $(i\, j\, k)$  of $(1\, 2\, 3)$ exists such that  $\alpha_{j}=0$ and $\alpha_{k}=0$.  In this case ${\bm{H}}_{\widetilde{n}}$ is normal to the face of $\widetilde{\mathcal{F}}$ belonging to the boundary of $K$ and corresponding to the face of $\widehat{K}$ of equation $\widehat{x}_i=1$. 
In order to be {$\mathrm{H}(\curl, {\mesh})$ conforming, ${\bm{H}}_{\widetilde{n}}$ is uniquely extended to all hexahedra $K'\in\mathcal{K}$ as follows:
if $K'$ contains the edge $\widetilde{E}\in\widetilde{\mathcal{E}}$ belonging to the boundary of $K$ and corresponding to the edge of $\widehat{K}$ of equations $\widehat{x}_{j}=1$, $\widehat{x}_{k}=1$, then  ${\bm{H}}_{\widetilde{n}}$ is extended to $K'$ by $\widehat{\bm{H}}_{{\widetilde{n}},\widehat{K}'}=\widehat{\mathbf{e}}_{i'}\widehat{\ell}_{\boldsymbol{\alpha}'}(-\widehat{\mathbf{x}})$, $\boldsymbol{\alpha}'$ being a permutation of $\boldsymbol{\alpha}$;
otherwise ${\bm{H}}_{\widetilde{n}}$ is extended by $\mathbf{0}$ to $K'$.
The number of  vector functions of this kind is $(P+1) |\widetilde{\mathcal{E}}|$.
\\

It is noticed that the support of each vector function $\bm{H}_{\widetilde{n}}$ in these three cases is contained in a single cell ${T}\in{\mesh}$. 

\subsubsection{{Semi-Discrete} Formulation}
A discrete formulation is provided by (\ref{E3}), (\ref{E4}) assuming that $\bm{E}$ and $\bm{e}$ belong to $H_P(\curl, \widetilde{\mesh})$ and that $\bm{H}$ and $\bm{h}$ belong
to $H_P(\curl, {\mesh})$. Thus assuming
\begin{align*}
& \bm{E}=\sum_{n=1}^N \bm{E}_n e_n\\
& \bm{H}=\sum_{\widetilde{n}=1}^{\widetilde{N}} \bm{H}_{\widetilde{n}} h_{\widetilde{n}},
\end{align*}
it results in
\begin{align}
& \mathbf{M}_{\varepsilon} \frac{d\mathbf{e}}{dt} = \widetilde{\mathbf{C}} \mathbf h,\label{F1}\\
  & \mathbf{M}_{\mu} \frac{d\mathbf{h}}{dt} = -{\mathbf{C}} \mathbf{e},\label{F2}
\end{align}
in which $\mathbf{e}$ is the column vector of the electric field degrees of freedom $e_n$, with $n=1,\dots,N$ and $\mathbf{h}$ is the column vector of the magnetic field degrees of freedom $h_{\widetilde{n}}$, with ${\widetilde{n}}=1,\dots,\widetilde{N}$.

The $N\times \widetilde{N}$ matrix $\widetilde{\mathbf{C}}$ has elements
\begin{align*}
	& \widetilde{c}_{n\widetilde{n}}=  \sum_{K \in \mathcal{K}} \left(-\int_{\partial K_{\dualmesh}} \bm{E}_n\times \bm{H}_{\widetilde{n}} \cdot \mathbf{n}\, d S+ 
 \int_{K} \curl\bm{E}_n \cdot \bm{H}_{\widetilde{n}}\, d V\right)
\end{align*}
with $n=1,\dots,N$ and $\widetilde{n}=1,\dots,\widetilde{N}$. The $ \widetilde{N}\times N$ matrix $\mathbf{C}$ has elements
\begin{align*}
	& {c}_{\widetilde{n}n}= \sum_{K \in \mathcal{K}} \left(-\int_{\partial K_{\mathcal{{T}}}} \bm{H}_{\widetilde{n}}\times \bm{E}_n \cdot \mathbf{n}\, d S +
 \int_{K} \curl\bm{H}_{\widetilde{n}} \cdot \bm{E}_n\, d V\right),
\end{align*}
with $\widetilde{n}=1,\dots,\widetilde{N}$ and $n=1,\dots,{N}$. By applying Stokes' theorem to this equation it thus follows $\widetilde{c}_{n\widetilde{n}}=c_{\widetilde{n}n}$ for  $n=1,\dots,{N}$ and $\widetilde{n}=1,\dots,\widetilde{N}$ or equivalently  $\widetilde{\mathbf{C}}=\mathbf{C}^T$.
It is noticed that element ${c}_{{n\widetilde{n}}}$ of $\mathbf{C}$ is determined by vector function $\bm{E}_n$ whose support is contained in one cell $\widetilde{T}\in\widetilde{\mesh}$ and by vector function $\bm{H}_{{\widetilde{n}}}$ whose support is contained in one cell ${T}\in{\mesh}$. Thus ${c}_{{n\widetilde{n}}}= 0$ only if $\dualtrig$ and $\trig$ are disjoint and matrix $\mathbf{C}$ is sparse.

Moreover from (\ref{E5}),  (\ref{E6}) it ensues
\begin{align*}
	& \widetilde{c}_{n\widetilde{n}}= \sum_{K \in \mathcal{K}} \left(-\int_{\partial \widehat{K}_{\dualmesh}} \widehat{\bm{E}}_{n,\widehat{K}} \times \widehat{\bm{H}}_{\widetilde{n},\widehat{K}} \cdot \widehat{\mathbf{n}}\, d\widehat{ S} + 
  \int_{\widehat{K}}\widehat{\curl}\widehat{\bm{E}}_{n,\widehat{K}}\cdot \widehat{\bm{H}}_{\widetilde{n},\widehat{K}}\, d\widehat{ V}\right), \\
  	& {c}_{\widetilde{n}n}= \sum_{K \in \mathcal{K}} \left(-\int_{\partial \widehat{K}_{\mesh}} \widehat{\bm{H}}_{\widetilde{n},K}\times \widehat{\bm{E}}_{n,\widehat{K}}\cdot \widehat{\mathbf{n}}\, d \widehat{ S} +
  \int_{\widehat{K}} \widehat{\curl}\widehat{\bm{H}}_{\widetilde{n},\widehat{K}}\cdot \widehat{\bm{E}}_{n,\widehat{K}}\, d\widehat{ V}\right),
\end{align*}
in which all integrals do not depend on the geometry or the materials. As a result all non-zero integrals for each $K$ can be computed only once for the single reference cube $\widehat{K}$.

The $N$-order  symmetric positive definite matrix $\mathbf{M}_{\varepsilon}$ has elements 
\begin{align}
M_{\varepsilon, mn} &=  \sum_{K \in \mathcal{K}} \int_{K} \bm{E}_m\cdot \boldsymbol{\varepsilon} \bm{E}_n\, d V =  \sum_{K \in \mathcal{K}} \int_{\widehat{K}} \widehat{\bm{E}}_{m,\widehat{K}} \cdot \widehat{\boldsymbol{\varepsilon}}_K\widehat{\bm{E}}_{n,\widehat{K}} \, |\mathbf{J}_{K}| d\widehat{ V}.\label{E9}
\end{align}
with $m,n=1,\dots,N$. By ordering the basis functions in such a way that the indexes $n$ of basis functions $\bm{E}_n$ whose support belong to the same cell  $\widetilde{T}\in\widetilde{\mesh}$ are consecutive,
$\mathbf{M}_{\varepsilon}$ is block diagonal.
The sparsity of $\mathbf{M}_{\varepsilon}$ is further increased by using, for the numerical integration of (\ref{E9}),  quadrature points $\widehat{\boldsymbol{\xi}}_{\boldsymbol{\alpha}}=(\widehat{\xi}_{\alpha_1}, \widehat{\xi}_{\alpha_2}, \widehat{\xi}_{\alpha_3})$ with weights $\widehat{W}_{\boldsymbol{\alpha}}=\widehat{w}_{\alpha_1}\cdot \widehat{w}_{\alpha_2}\cdot \widehat{w}_{\alpha_3}$, being $\boldsymbol{\alpha}=(\alpha_1, \alpha_2, \alpha_3)\in\{0,\dots,P\}^3$.
In fact by assuming that in the reference cube $\mathcal{K}$
 it is $\widehat{\bm{E}}_{m,\widehat{K}}=\widehat{\ell}_{\boldsymbol{\beta}}(\widehat{\mathbf{x}})\widehat{\mathbf{e}}_j$
and $\widehat{\bm{E}}_{n,\widehat{K}}=\widehat{\ell}_{\boldsymbol{\alpha}}(\widehat{\mathbf{x}})\widehat{\mathbf{e}}_i$ the integral 
\begin{align*}
& \int_{\widehat{K}} \widehat{\bm{E}}_{m,\widehat{K}} \cdot \widehat{\boldsymbol{\varepsilon}}_K\widehat{\bm{E}}_{n,\widehat{K}} \, |\mathbf{J}_{K}| d\widehat{ V}. 
\end{align*}
is approximated by
\begin{align}
&  \widehat{W}_{\boldsymbol{\alpha}} |\mathbf{J}_K(\widehat{\boldsymbol{\xi}}_{\boldsymbol{\alpha}})|\widehat{\boldsymbol{\varepsilon}}_{K,ij}(\widehat{\boldsymbol{\xi}}_{\boldsymbol{\alpha}}) \,\delta_{\boldsymbol{\alpha\beta}}. \label{E10}
\end{align}
$\delta$ being Kronecker's symbol.
Hence in all hexahedra $K'\in\mathcal{K}$ in which both $\widehat{\bm{E}}_{m}$ and  $\widehat{\bm{E}}_{n}$ are not identically zero and thus given by $\widehat{\bm{E}}_{m,\widehat{K}'}=\widehat{\ell}_{\boldsymbol{\beta}'}(\widehat{\mathbf{x}})\widehat{\mathbf{e}}_{j'}$ and $\widehat{\bm{E}}_{n,\widehat{K}'}=\widehat{\ell}_{\boldsymbol{\alpha}'}(\widehat{\mathbf{x}})\widehat{\mathbf{e}}_{i'}$,  $\boldsymbol{\alpha}'$ and $\boldsymbol{\beta}'$ are the same permutation of $\boldsymbol{\alpha}$ and $\boldsymbol{\beta}$ respectively. As a result if $\boldsymbol{\alpha}\neq\boldsymbol{\beta}$ in $K$ then $\boldsymbol{\alpha}'\neq\boldsymbol{\beta}'$ in each $K'$ and, from  (\ref{E9}) and (\ref{E10}),
 $M_{\varepsilon, mn}=0$.

{In the following we refer to $\mathbf{M}_\varepsilon$ and $\mathbf{M}_\mu$ defined in \eqref{E9} and \eqref{E12} as the \emph{consistent mass matrices}, since they arise from the exact (or highly accurate) numerical integration of the $L^2$ inner products on primal and dual cells. The Gauss--Radau-based tensor-product quadrature rule introduced in \eqref{E10}--\eqref{E13} is exact for polynomials up to degree $2P$ in each coordinate direction, which is sufficient to preserve the nominal spatial convergence order of the method. When we additionally exploit the Kronecker structure of the basis and quadrature to discard off-diagonal couplings between different interpolation nodes, we obtain \emph{lumped mass matrices}, which are block-diagonal with small dense blocks corresponding to the supports of basis functions. Unless stated otherwise, the term ``mass matrices'' refers generically to either the consistent or lumped versions, while ``lumped mass matrices'' always denotes the block-diagonal approximations used for explicit time stepping and for the direct application of $\mathbf{M}_\varepsilon^{-1}$ and $\mathbf{M}_\mu^{-1}$.}

The $\widetilde{N}$-order  symmetric positive definite matrix $\mathbf{M}_{\mu}$ has elements 
\begin{align}\label{E12}
M_{\mu, \widetilde{m}\widetilde{n}} &=  \sum_{K \in \mathcal{K}} \int_{K} \bm{H}_{\widetilde{m}}\cdot \boldsymbol{\mu} \bm{H}_{\widetilde{n}}\, d V =  \sum_{K \in \mathcal{K}} \int_{\widehat{K}} \widehat{\bm{H}}_{{\widetilde{m}},\widehat{K}} \cdot \widehat{\boldsymbol{\mu}}_K\widehat{\bm{H}}_{{\widetilde{n}},\widehat{K}} \, |\mathbf{J}_{K}| d\widehat{ V}.
\end{align}
with ${\widetilde{m}},{\widetilde{n}}=1,\dots,\widetilde{N}$. By ordering the basis functions in such a way that the indexes ${\widetilde{n}}$ of basis functions $\bm{H}_{\widetilde{n}}$ whose support belong to the same  cell  ${T}\in{\mesh}$ are consecutive
$\mathbf{M}_{\mu}$ is block diagonal.
The sparsity of $\mathbf{M}_{\mu}$ if further increased by using, for the numerical integration of (\ref{E12}),  quadrature points $-\widehat{\boldsymbol{\xi}}_{\boldsymbol{\alpha}}$ with weights $\widehat{Z}_{\boldsymbol{\alpha}}=\widehat{W}_{\boldsymbol{\alpha}}$, being $\boldsymbol{\alpha}\in\{0,\dots,P\}^3$.
In fact by assuming that in the reference cube $\widehat{K}$ it is $\widehat{\bm{H}}_{{\widetilde{m}},\widehat{K}}=\widehat{\ell}_{\boldsymbol{\beta}}(-\widehat{\mathbf{x}})\widehat{\mathbf{e}}_j$ and $\widehat{\bm{H}}_{{\widetilde{n}},\widehat{K}}=\widehat{\ell}_{\boldsymbol{\alpha}}(-\widehat{\mathbf{x}})\widehat{\mathbf{e}}_i$ the integral
\begin{align*}
& \int_{\widehat{K}} \widehat{\bm{H}}_{\widetilde{m},\widehat{K}} \cdot \widehat{\boldsymbol{\mu}}_K\widehat{\bm{H}}_{\widetilde{n},\widehat{K}} \, |\mathbf{J}_{K}| d\widehat{ V} 
\end{align*}
can be approximated by
\begin{align}
& \widehat{Z}_{\boldsymbol{\alpha}} |\mathbf{J}_K(-\widehat{\boldsymbol{\xi}}_{\boldsymbol{\alpha}})|\widehat{\boldsymbol{\mu}}_{K,ij}(-\widehat{\boldsymbol{\xi}}_{\boldsymbol{\alpha}}) \,\delta_{\boldsymbol{\alpha\beta}}. \label{E13}
\end{align}
Hence in all hexahedra $K'\in\mathcal{K}$ in which both ${\bm{H}}_{\widetilde{m}}$ and  ${\bm{H}}_{\widetilde{n}}$ are not identically zero so that $\widehat{\bm{H}}_{\widetilde{m},\widehat{K}'}=\widehat{\ell}_{\boldsymbol{\beta}'}(-\widehat{\mathbf{x}})\widehat{\mathbf{e}}_{j'}$ and $\widehat{\bm{H}}_{\widetilde{n},\widehat{K}'}=\widehat{\ell}_{\boldsymbol{\alpha}'}(-\widehat{\mathbf{x}})\widehat{\mathbf{e}}_{i'}$,  $\boldsymbol{\alpha}'$ and $\boldsymbol{\beta}'$ are the same permutation of $\boldsymbol{\alpha}$ and $\boldsymbol{\beta}$ respectively. As a result if $\boldsymbol{\alpha}\neq\boldsymbol{\beta}$ in $K$ then $\boldsymbol{\alpha}'\neq\boldsymbol{\beta}'$ in each $K'$, and from  (\ref{E12}) and (\ref{E13}), $M_{\mu, \widetilde{m}\widetilde{n}}=0$.

%
\subsection{Time Discretisation}\label{sec:timedisc}

{For the temporal discretisation of the semi-discrete system \eqref{F1}--\eqref{F2} we employ the classical leap-frog scheme. This choice is natural in combination with the skew-adjoint structure of the discrete curl operator and the block-diagonal mass matrices: the method is explicit, second-order accurate in time, and preserves a discrete analogue of the electromagnetic energy in the lossless case.}
The common leap frog scheme for the time discretisation of (\ref{F1}), (\ref{F2}) is used as a first choice:
\begin{align}
  \begin{split}
  & \mathbf{h}^{1/2} =\mathbf{h}^{0} -\frac{1}{2}\Delta t \,\mathbf{M}_{\boldsymbol{\mu}}^{-1} \mathbf{C} \mathbf{e}^0,\\
& \mathbf{e}^{q+1} = \mathbf{e}^{q} +\Delta t \,\mathbf{M}_{\boldsymbol{\varepsilon}}^{-1} \mathbf{C}^T \mathbf{h}^{q+1/2},\\
& \mathbf{h}^{q+1/2} =\mathbf{h}^{q-1/2} -\Delta t \,\mathbf{M}_{\boldsymbol{\mu}}^{-1} \mathbf{C} \mathbf{e}^q,
  \end{split}
  \label{eq:leap_frog}
\end{align}
in which $\mathbf{e}^q$ is the approximation of $\mathbf{e}$ at time instant $q\,\Delta t$ and $\mathbf{h}^{q+1/2}$ is the approximation of $\mathbf{h}$  at time instant $(q+1/2)\,\Delta t$, for $q=0,\dots,Q$.

In order to guarantee stability, the time step $\Delta t$ is computed in such a way that $\Delta t < 2 / \sqrt{\lambda_M}$, being $\lambda_M$ the maximum of the eigenvalues $\lambda$ of the eigenvalue problem
\begin{align}
& \mathbf{C} \mathbf{M}_{\boldsymbol{\varepsilon}}^{-1} \mathbf{C}^T \mathbf{u} = \lambda \mathbf{M}_{\boldsymbol{\mu}} \mathbf{u}
  \label{eq:evp}
\end{align} 
for the eigenvectors $\mathbf{u}$.

\section{Numerical examples}\label{sec:numerics}

{In this section we present a series of numerical experiments that illustrate the accuracy, spectral properties, and computational efficiency of the proposed three-dimensional dual cell method. We first investigate the spectrum of the discrete curl--curl operator, demonstrating spectral correctness (absence of spurious eigenvalues) and high-order convergence of selected eigenvalues and eigenfunctions. We then assess the stability constraint imposed by the CFL condition and study the convergence of time-domain simulations in straight and curved waveguides, including a configuration with a spherical dielectric inclusion and perfectly matched layers (PMLs)~\cite{collinoOptimizingPerfectlyMatched1998,teixeiraComplexSpaceApproach2000}. Finally, we analyse the sparsity pattern of the (lumped) mass matrices and verify that the number of non-zero entries per row remains essentially independent of the polynomial order, thereby confirming the expected efficiency gains of the lumped formulation.}

{All experiments are performed using an implementation of the method in the open-source high-order finite element library Netgen/NGsolve~\cite{ngsolve}. The dual cell method extension is available as an add-on package~\cite{dcm}. Discrete eigenvalue problems are solved using the time-domain Krylov approach from~\cite{NannenWess2024}, which, in the present setting, only requires repeated applications of the stiffness matrix and of the inverse mass matrices. Owing to the block-diagonal structure of $\mathbf{M}_\varepsilon$ and $\mathbf{M}_\mu$ produced by our lumping strategy, these inverses are applied by direct inversion of small local blocks rather than by an additional Krylov iteration, in contrast to standard FEM discretisations with globally coupled mass matrices~\cite{jinFiniteElementMethod2014a,weilandTimeDomainElectromagnetic1996}.}

\subsection{Sparsity of the mass matrices}

\begin{figure}
  \centering
  \includegraphics{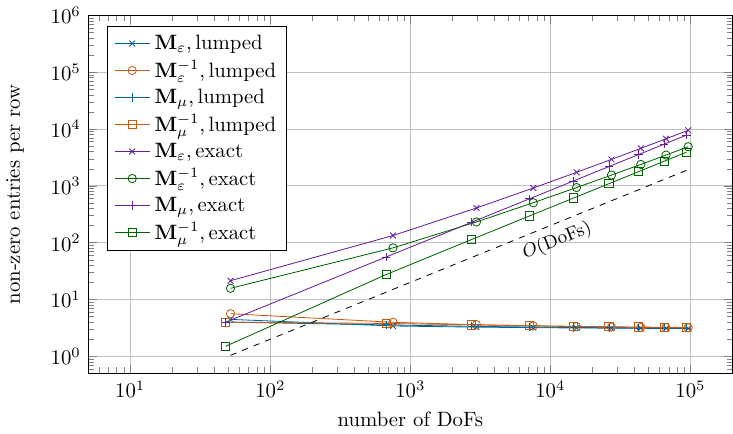}
  \caption{Sparsity of the resulting mass matrices for a fixed mesh and increasing polynomial order.}
  \label{fig:sparsity}
\end{figure}
{The goal of the employed mass-lumping strategy is to obtain lumped mass matrices whose sparsity pattern does not deteriorate with increasing polynomial order. To study this behaviour we consider the unit cube $\Omega = (0,1)^3$ discretised by four tetrahedra and construct the associated barycentric dual mesh. For a sequence of polynomial orders $P$ we assemble both the consistent mass matrices $\mathbf{M}_\varepsilon$, $\mathbf{M}_\mu$ (as in \eqref{E9}--\eqref{E13}) and their lumped counterparts, and we record the number of non-zero entries per row. Figure~\ref{fig:sparsity} shows that, for the lumped mass matrices, the number of non-zero entries per row approaches a value close to three as $P$ increases, reflecting the locality and near-orthogonality of the Gauss--Radau-based basis functions within each dual cell. In contrast, the number of non-zero entries per row in the consistent mass matrices grows essentially linearly with the total number of global degrees of freedom, as expected from a standard high-order edge-element discretisation.}

\subsection{Eigenvalue problem}
\begin{figure}
  \begin{subfigure}[t]{0.5\textwidth}
    \includegraphics[width=\textwidth]{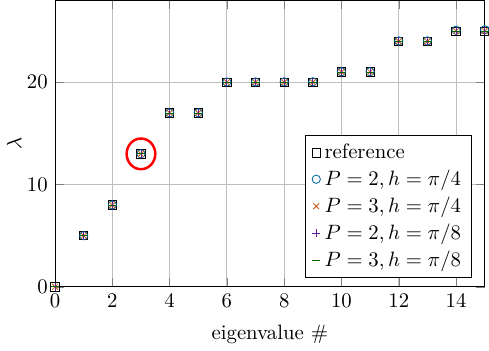}
    \caption{Spectrum for selected discretisations}
    \label{fig:evs}
  \end{subfigure}
  \begin{subfigure}[t]{0.5\textwidth}
     \includegraphics[width = \textwidth]{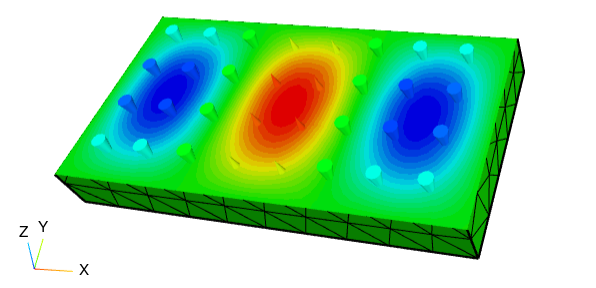}
     \caption{Approximated eigenfunction corresponding to the second non-trivial eigenvalue $\lambda_3 = 13$ of \eqref{eq:evp} on the domain $\Omega$ from \eqref{eq:evp_domain} on a typical mesh.}
  \label{fig:ef}
  \end{subfigure}
\end{figure}

\begin{figure}
  \begin{subfigure}[t]{0.5\textwidth}
  \includegraphics[width=\textwidth]{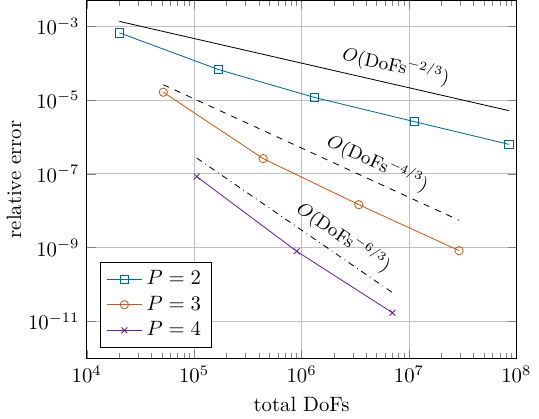}
  \caption{Convergence with respect to the total number of DoFs.}
  \label{fig:ev_conv_dofs_3}
  \end{subfigure}
  \begin{subfigure}[t]{0.5\textwidth}
  \includegraphics[width=\textwidth]{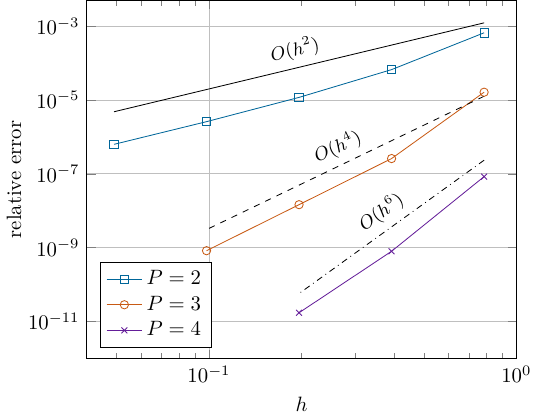}
  \caption{Convergence with respect to the mesh size $h$.}
  \label{fig:ev_conv_h_3}
  \end{subfigure}
  \caption{Convergence of the third, non-trivial eigenvalue.}
  \label{fig:ev_conv_3}
\end{figure}

\begin{figure}
  \begin{subfigure}[t]{0.5\textwidth}
  \includegraphics[width=\textwidth]{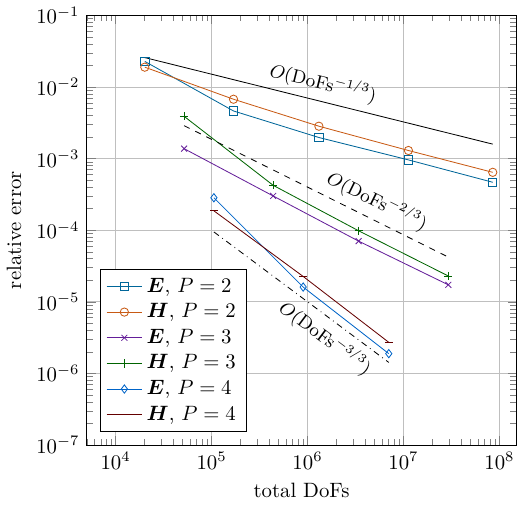}
  \caption{Convergence with respect to the total number of DoFs.}
  \label{fig:ef_conv_dofs_3}
  \end{subfigure}
  \begin{subfigure}[t]{0.5\textwidth}
  \includegraphics[width=\textwidth]{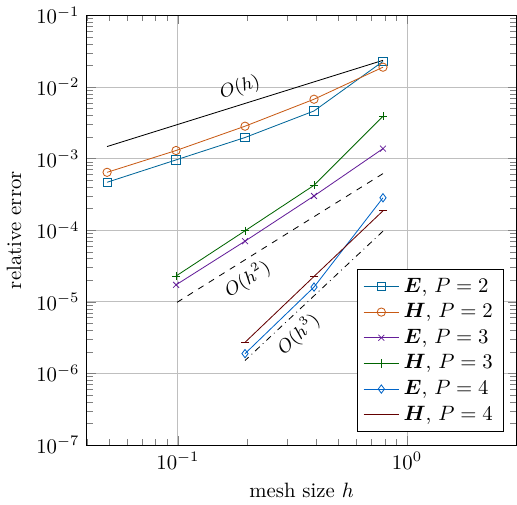}
  \caption{Convergence with respect to the mesh size $h$.}
  \label{fig:ef_conv_h_3}
  \end{subfigure}
  \caption{Convergence of the eigenfunctions corresponding to the second, non-trivial eigenvalue}
  \label{fig:ef_conv_3}
\end{figure}

{We next analyse the discrete eigenvalue problem \eqref{eq:evp} on the rectangular cavity (see also Figure~\ref{fig:ef})}
\begin{align}
  \Omega = (0,\pi)\times(0,\pi/2)\times(0,\pi/4),
  \label{eq:evp_domain}
\end{align}
  {with perfectly conducting boundary conditions $\bm{E}\times \bm{n} = \bm{0}$ on $\partial \Omega$.}

  {Figure~\ref{fig:evs} shows the resulting discrete spectra for selected discretisations with different mesh sizes $h$ and polynomial orders $P$. All computed discrete eigenvalues approximate analytically known eigenvalues of the continuous cavity problem, and no spurious eigenvalues are observed; this corroborates the spectral correctness of the three-dimensional extension, in agreement with the two-dimensional results of~\cite{kapidaniArbitraryorderCellMethod2021}.}

  {Figure~\ref{fig:ev_conv_3} reports the relative error of the discrete eigenvalues approximating the second non-trivial eigenvalue $\lambda_2=8$ with respect to the total number of unknowns for several polynomial orders $P$. The corresponding eigenmode is depicted in Figure~\ref{fig:ef}. We observe convergence rates of order $O\!\left(\mathrm{DoFs}^{(2P-2)/3}\right)$, which translates to $O\!\left(h^{2P-2}\right)$ with respect to the mesh size $h$, consistent with the $(P-1)$-degree approximation of the curl-conforming fields.}

  {In Figure~\ref{fig:ef_conv_3} we show the convergence of the relative errors in the (mass-lumped) $L_2$-norms of the electric and magnetic fields of the eigenfunctions corresponding to $\lambda_2=8$. The errors are measured in the mass-lumped $L^2$ norms induced by $\mathbf{M}_\varepsilon$ and $\mathbf{M}_\mu$, respectively. For both fields we observe convergence of order $O\!\left(\mathrm{DoFs}^{(P-1)/3}\right)$, i.e., $O\!\left(h^{P-1}\right)$ with respect to $h$, in agreement with the theoretical approximation properties of the underlying spaces.}

\subsection{The time-domain problem}
In this section we present numerical experiments, examining the fully discretised time-domain problem \eqref{eq:leap_frog}.
\subsubsection{CFL condition}
\begin{figure}
  \begin{subfigure}[t]{0.5\textwidth}
      \includegraphics[width=\textwidth]{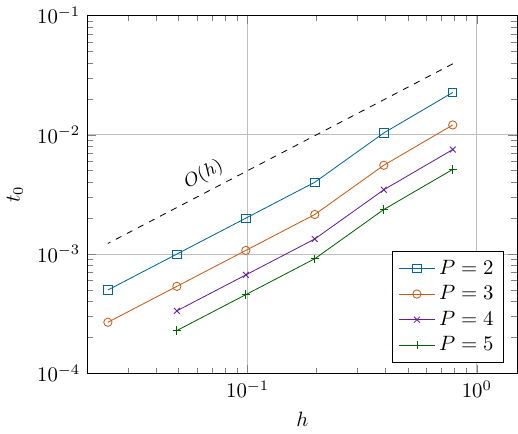}
      \caption{The maximal stable time-step for fixed values of $P$ and varying mesh-size $h$.}
  \end{subfigure}
  \begin{subfigure}[t]{0.5\textwidth}
      \includegraphics[width=\textwidth]{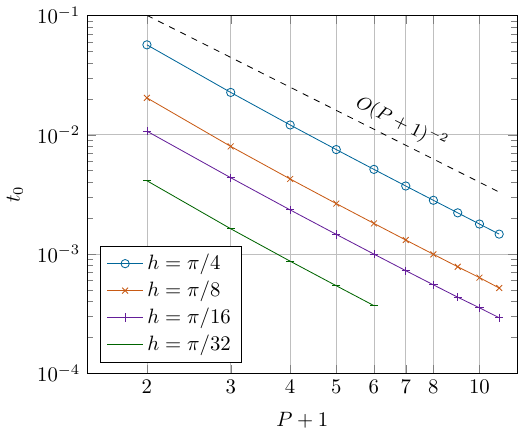}
      \caption{The maximal stable time-step for fixed mesh-sizes $h$ and varying polynomial order $P$.}
  \end{subfigure}
  \caption{The CFL condition with respect to mesh-size and polynomial order.}
  \label{fig:cfl}
\end{figure}

{We first verify the dependence of the CFL stability bound on the mesh size $h$ and on the polynomial order $P$ using the cavity problem from the previous subsection. For each pair $(h,P)$ we approximate the largest eigenvalue $\lambda_M$ of the generalised eigenvalue problem~\eqref{eq:evp} by a simple power iteration applied to the matrix pencil $(\mathbf{C}\mathbf{M}_{\boldsymbol{\varepsilon}}^{-1}\mathbf{C}^T,\mathbf{M}_{\boldsymbol{\mu}})$ and compute the corresponding maximal allowed time step $t_0 \approx 2/\sqrt{\lambda_M}$. Figure~\ref{fig:cfl} displays the resulting stable time steps. We observe that $t_0$ scales linearly with $h$ and decreases like $(P+1)^{-2}$, in agreement with standard dispersion analyses for high-order explicit time-domain discretisations~\cite{hesthavenNodalDiscontinuousGalerkin2008}.}

\subsubsection{Convergence of the time-domain solution}
\label{sec:waveguide}
\begin{figure}
  \begin{subfigure}[t]{0.5\textwidth}
    \includegraphics[width=\textwidth, clip, trim = 0 0 0 100]{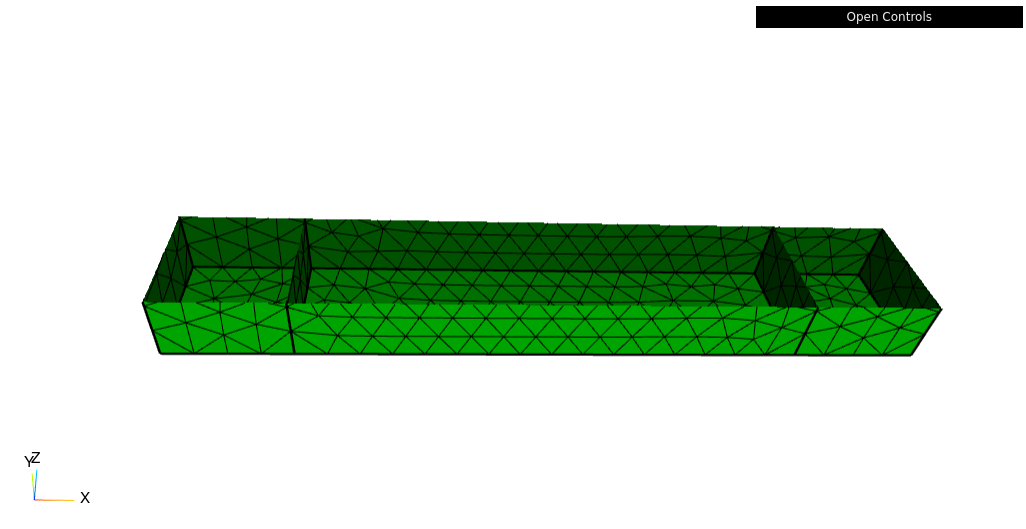}
      \caption{Empty waveguide (Subsection \ref{sec:waveguide}).}
  \label{fig:waveguide_mesh}
  \end{subfigure}
  \begin{subfigure}[t]{0.5\textwidth}
    \includegraphics[width=\textwidth, clip, trim = 0 0 0 100]{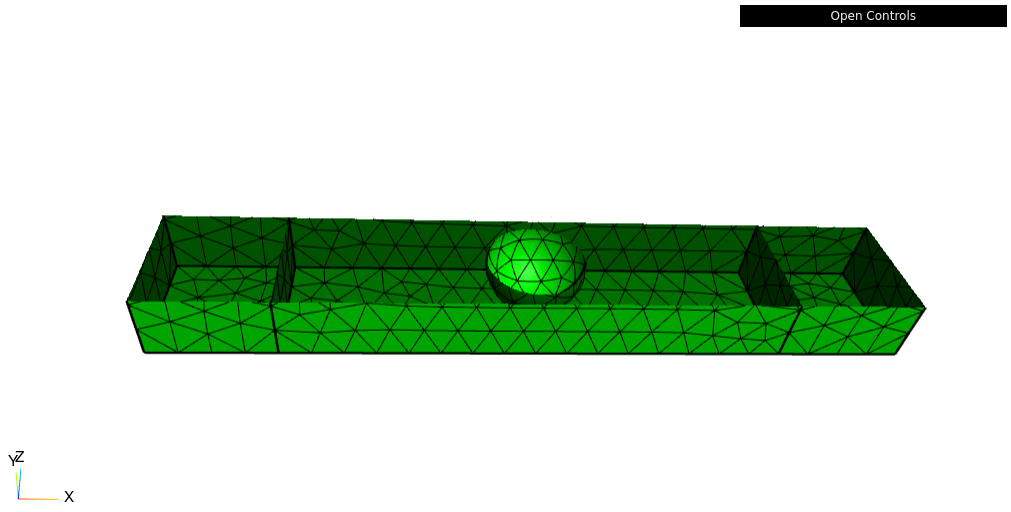}
      \caption{Waveguide with spherical inclusion (Subsection \ref{sec:waveguide_sphere}).}
  \label{fig:waveguide_sphere_mesh}
  \end{subfigure}
  \caption{Typical meshes of the waveguide problems with perfectly matched layers.}
  \label{fig:waveguide_meshes}
\end{figure}
To show time-domain convergence we use a waveguide (cf. Figure \ref{fig:waveguide_mesh})
\begin{align*}
  \Omega_{int}=(0,2)\times(0,1/2)\times(0,1/2),
\end{align*}
with attached perfectly matched layers left and right with length $l=0.5$ (for an exact formulation of the PML problem see the Appendix). We pick constant material parameters $\varepsilon\equiv\mu\equiv 1$. We propagate a wave 
\begin{align}
  E_0(t,y,z)=\left(0,0,\exp(-5(1-t)^2)\sin(10t)\sin(2\pi y)\right),
  \label{eq:incoming_wave}
\end{align}
from the left side $\Gamma_{in}=\{0\}\times(0,1/2)\times(0,1/2)$ of the waveguide. 
{In all time-domain examples we choose the time-step size $\Delta t$ according to the CFL bound of the finest discretisation under consideration. This guarantees stability for all coarser meshes and ensures that the temporal error remains subdominant compared to the spatial discretisation error.}

{Figure~\ref{fig:waveguide} shows snapshots of the electric field magnitude in the interior waveguide region at selected times. A semi-analytic reference solution for this configuration can be computed by Laplace-domain techniques (see Appendix~\ref{sec:ref_sol}); the resulting convolution in time is evaluated numerically. In Figure~\ref{fig:td_waveguide_error} we report the $L_2(\Omega_{int})$-error of the electric field over time for polynomial order $P=2$ and several mesh sizes $h$. The plots confirm a monotone decrease of the error as the mesh is refined and show that no additional oscillations or instabilities are introduced by the explicit time stepping.}
\begin{figure}
  \begin{subfigure}{0.5\textwidth}
      \includegraphics[width=\textwidth]{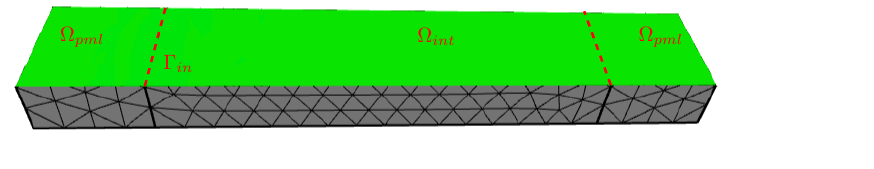}
      \caption{$t=0.5$}
  \end{subfigure}
  \begin{subfigure}{0.5\textwidth}
      \includegraphics[width=\textwidth]{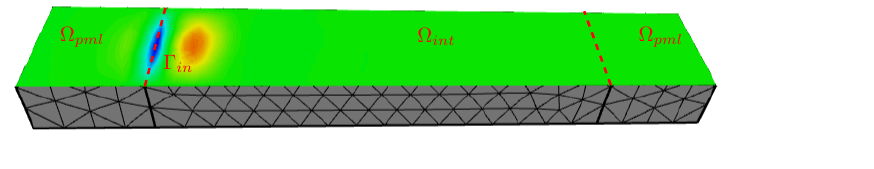}
      \caption{$t=1.0$}
  \end{subfigure}
  \begin{subfigure}{0.5\textwidth}
      \includegraphics[width=\textwidth]{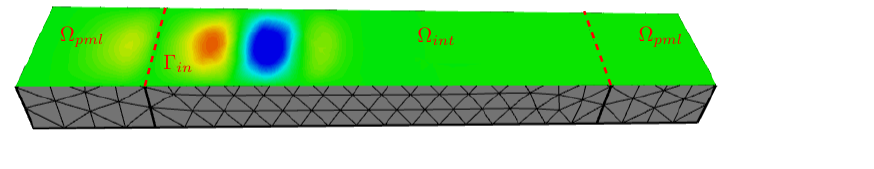}
      \caption{$t=1.5$}
  \end{subfigure}
  \begin{subfigure}{0.5\textwidth}
      \includegraphics[width=\textwidth]{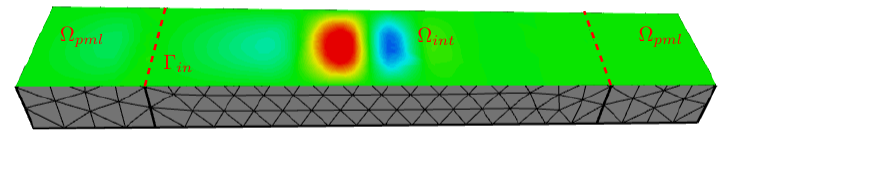}
      \caption{$t=2.0$}
  \end{subfigure}
  \begin{subfigure}{0.5\textwidth}
      \includegraphics[width=\textwidth]{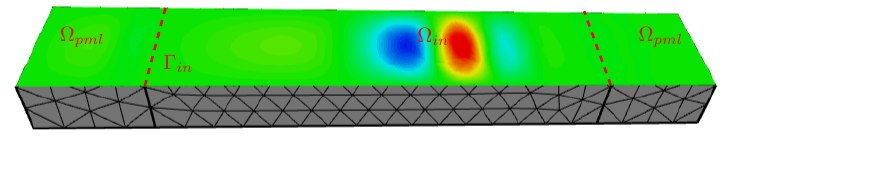}
      \caption{$t=2.5$}
  \end{subfigure}
  \begin{subfigure}{0.5\textwidth}
      \includegraphics[width=\textwidth]{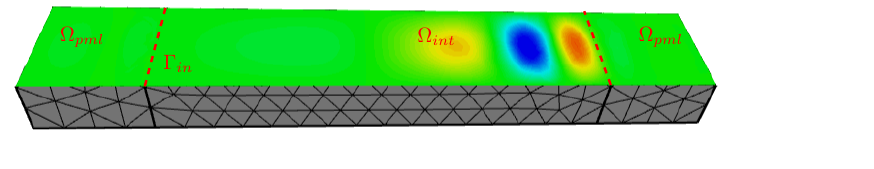}
      \caption{$t=3.0$}
  \end{subfigure}
  \begin{subfigure}{0.5\textwidth}
      \includegraphics[width=\textwidth]{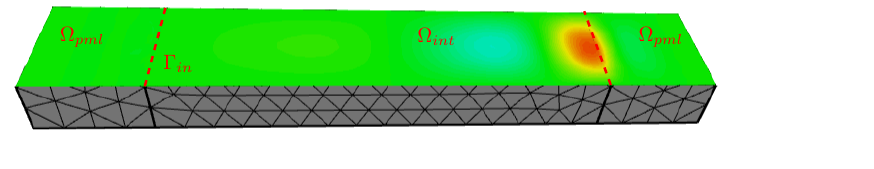}
      \caption{$t=3.5$}
  \end{subfigure}
  \begin{subfigure}{0.5\textwidth}
      \includegraphics[width=\textwidth]{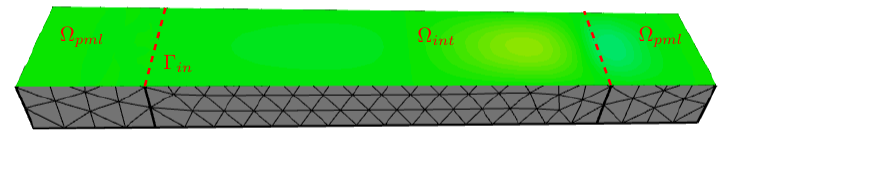}
      \caption{$t=4.0$}
  \end{subfigure}
  \caption{Snapshots of the plain waveguide problem.}
  \label{fig:waveguide}
\end{figure}

\begin{figure}
  \begin{subfigure}{0.5\textwidth}
      \includegraphics[width=\textwidth]{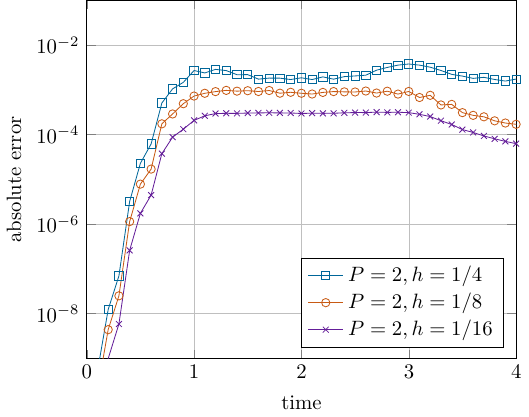}
      \caption{$P=2$}
  \end{subfigure}
  \begin{subfigure}{0.5\textwidth}
      \includegraphics[width=\textwidth]{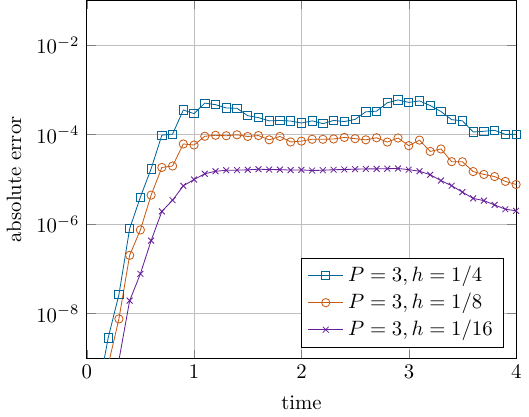}
      \caption{$P=3$}
  \end{subfigure}
  \caption{$L_2(\Omega_{int})$-error of the $E$ field (cf. Figure \ref{fig:waveguide}) over time for fixed polynomial order and different mesh-sizes $h$.}
  \label{fig:td_waveguide_error}
\end{figure}

{Figure~\ref{fig:td_waveguide_conv} shows the convergence of the relative space--time $L_2([0,1];\Omega_{int})$-error of the electric field for various orders $P$. In accordance with the eigenfunction convergence results, we observe rates of order $O\!\left(\mathrm{DoFs}^{(P-1)/3}\right)$, i.e., $O\!\left(h^{P-1}\right)$ with respect to the mesh size. This indicates that, for the chosen $\Delta t$, the overall error is dominated by the spatial approximation.}

\begin{figure}
  \begin{subfigure}{0.5\textwidth}
      \includegraphics[width=\textwidth]{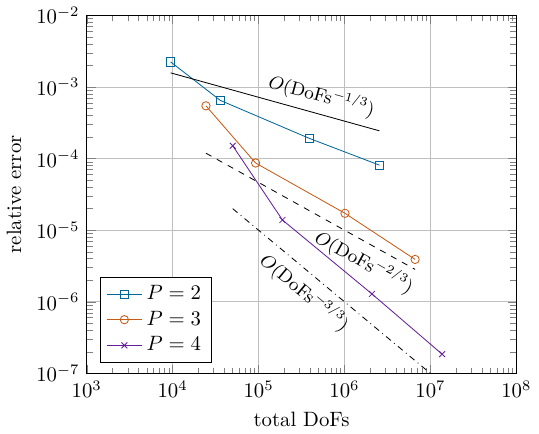}
      \caption{Convergence with respect to the total number of DoFs}
  \end{subfigure}
  \begin{subfigure}{0.5\textwidth}
      \includegraphics[width=\textwidth]{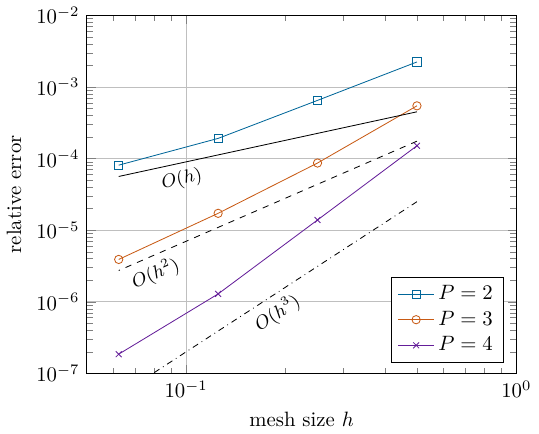}
      \caption{Convergence with respect to the mesh size $h$}
  \end{subfigure}
  \caption{Convergence of the $L_2([0,1];\Omega_{int})$-error of the time domain $E$ field (cf. Figure \ref{fig:waveguide}) for different polynomial orders $P$.}
  \label{fig:td_waveguide_conv}
\end{figure}

\subsubsection{Time domain problem in curved geometry}
\label{sec:waveguide_sphere}
To show that our method is also applicable to non-constant material parameters and curved geometries, we enhance the waveguide problem above by adding a dielectric sphere with $\varepsilon\equiv 9$ in the center of the waveguide (cf. Figure \ref{fig:waveguide_sphere_mesh}). Outside of the sphere we pick $\varepsilon\equiv 1$ and $\mu\equiv 1$ everywhere. {Figure~\ref{fig:waveguide_sphere} shows snapshots of the time-domain solution and illustrates the interaction of the guided wave with the dielectric inclusion and the subsequent scattering pattern.}

\begin{figure}
  \begin{subfigure}{0.5\textwidth}
      \includegraphics[width=\textwidth]{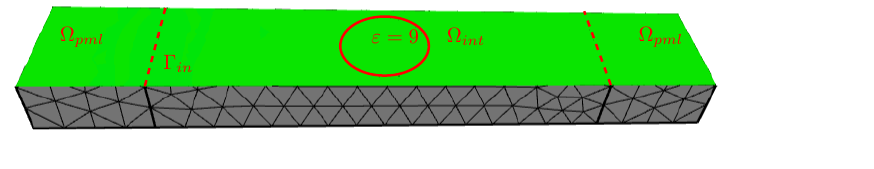}
      \caption{$t=0.5$}
  \end{subfigure}
  \begin{subfigure}{0.5\textwidth}
      \includegraphics[width=\textwidth]{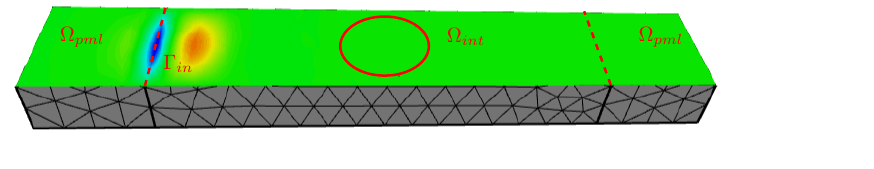}
      \caption{$t=1.0$}
  \end{subfigure}
  \begin{subfigure}{0.5\textwidth}
      \includegraphics[width=\textwidth]{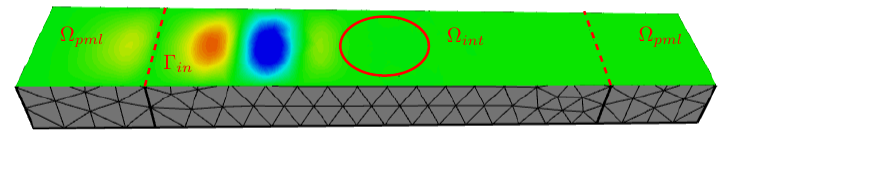}
      \caption{$t=1.5$}
  \end{subfigure}
  \begin{subfigure}{0.5\textwidth}
      \includegraphics[width=\textwidth]{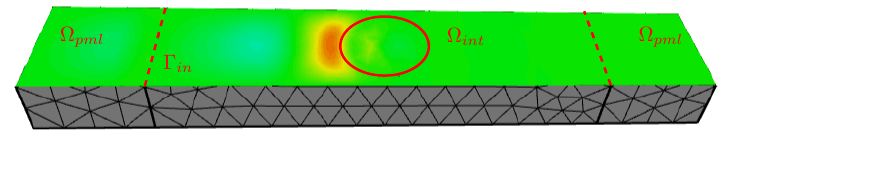}
      \caption{$t=2.0$}
  \end{subfigure}
  \begin{subfigure}{0.5\textwidth}
      \includegraphics[width=\textwidth]{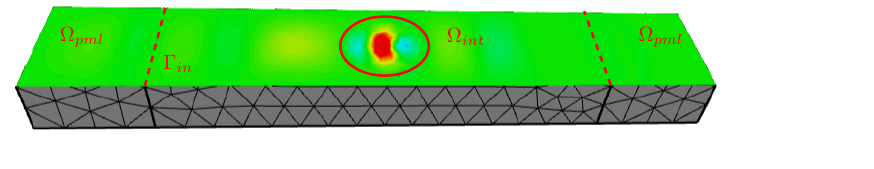}
      \caption{$t=2.5$}
  \end{subfigure}
  \begin{subfigure}{0.5\textwidth}
      \includegraphics[width=\textwidth]{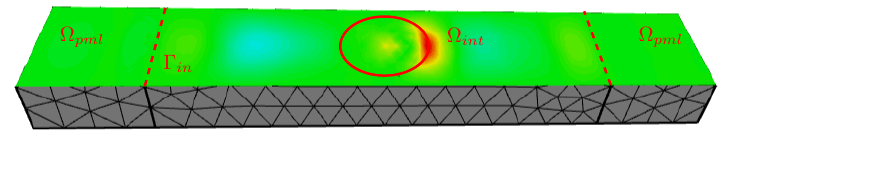}
      \caption{$t=3.0$}
  \end{subfigure}
  \begin{subfigure}{0.5\textwidth}
      \includegraphics[width=\textwidth]{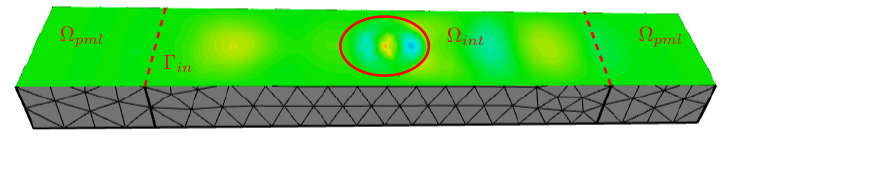}
      \caption{$t=3.5$}
  \end{subfigure}
  \begin{subfigure}{0.5\textwidth}
      \includegraphics[width=\textwidth]{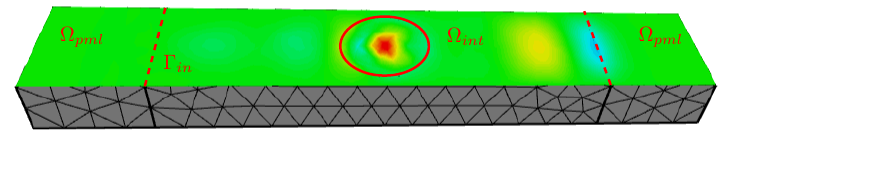}
      \caption{$t=4.0$}
  \end{subfigure}
  \begin{subfigure}{0.5\textwidth}
      \includegraphics[width=\textwidth]{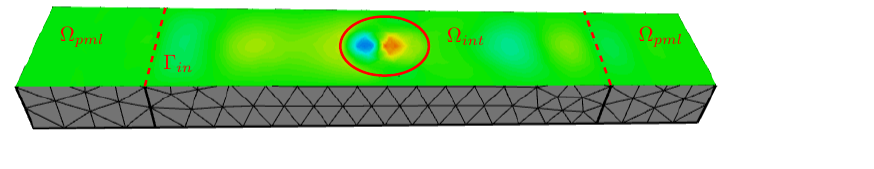}
      \caption{$t=4.5$}
  \end{subfigure}
  \begin{subfigure}{0.5\textwidth}
      \includegraphics[width=\textwidth]{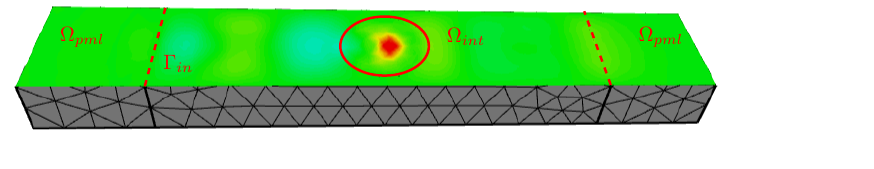}
      \caption{$t=5.0$}
  \end{subfigure}
  \caption{Waveguide with spherical inclusion}
  \label{fig:waveguide_sphere}
\end{figure}

\section{Conclusion}

{In this work we have extended the high-order dual cell method for time-domain Maxwell equations from two to three spatial dimensions. The construction relies on barycentric dual meshes and on pairs of mutually dual curl-conforming spaces, defined via tensor-product Gauss--Radau interpolation polynomials and mapped from a reference cube by trilinear covariant Piola transformations within each hexahedral subcell of primal and dual mesh elements. A key feature of the approach is the block-diagonal structure of the mass matrices, which allows for explicit time integration by a leap-frog scheme while preserving a discrete electromagnetic energy identity.}

{Through a series of numerical experiments we have demonstrated that the three-dimensional dual cell method inherits the favourable properties previously observed in two dimensions. In particular, the sparsity pattern of the lumped mass matrices remains essentially independent of the polynomial order, the discrete spectrum of the curl--curl operator is free of spurious modes and exhibits high-order convergence of eigenvalues and eigenfunctions, and the observed CFL bound scales linearly with the mesh size and quadratically with $(P+1)^{-1}$. Time-domain simulations in straight and curved waveguides, including configurations with perfectly matched layers and a spherical dielectric inclusion, confirm high-order convergence of the fields and show that complex geometries and material interfaces can be treated accurately and robustly.}

{Several extensions of the present work are natural directions for future research. An important next step is the treatment of more general boundary conditions and source terms, as well as dispersive or non-linear material laws, within the same dual cell framework. Further topics of interest include adaptive mesh refinement, local time stepping on barycentric dual meshes, and large-scale parallel implementations targeting multi-GPU and exascale architectures. Finally, coupling the dual cell method to circuit or lumped-element models, for instance in the context of microwave devices, appears to be a promising application area.}

\section*{Appendix}
\subsection{PML formulation}
\label{sec:pml}
In the numerical experiments in Sections \ref{sec:waveguide}, and \ref{sec:waveguide_sphere} we use a PML in a three-dimensional waveguide of the following weak formulation. To the interior waveguide $\Omega_{int}:=(0,a)\times(0,b)\times(0,c)$ we append layers $\Omega_{pml}:=\left((-l,0)\cup(a,a+l)\right)\times(0,b)\times(0,c)$.
The usual derivation of time-domain PMLs is done in frequency domain, where a mapping depending on the frequency is introduced. We avoid this detour here and nonchalantly write the frequency-dependency in time-domain as a dependency on $\partial_t$. 
Using the mapping 
\begin{align*}
  \tilde{\mathbf x}(\mathbf x)= \tilde x(x_1,x_2,x_3)=\mathbf x+\begin{cases}
    \frac{\sigma}{\partial_t}(x_1-a)e_1,&x_1>a,\\
    \frac{\sigma}{\partial_t}x_1e_1,&x_1<0,\\
    0,&\text{else,}
  \end{cases}
\end{align*}
and according Jacobian (inverse) matrix and determinant 
\begin{align*}
  \mathbf J(x) &= \mathbf I+
  \begin{cases}
    \frac{\sigma}{\partial_t}\Pi_1,&\mathbf x\in\Omega_{pml},\\
    0,&\text{else,}
  \end{cases},&
  \mathbf J(x)^{-1} &= 
  \begin{cases}
    (1+\frac{\sigma}{\partial_t})^{-1}\Pi_1+\Pi_1^\perp,&\mathbf x\in\Omega_{pml},\\
    \mathbf I,&\text{else,}
  \end{cases}\\
  \det \mathbf J(x) &= 1+
  \begin{cases}
    \frac{\sigma}{\partial_t},&\mathbf x\in\Omega_{pml},\\
    0,&\text{else,}
  \end{cases}
\end{align*}
where $\Pi_1 = \mathbf e_1\mathbf e_1^T$ is the orthogonal projection onto the $x$-component and $\Pi_1^\perp = \mathbf I-\Pi_1$.
Introducing $\bm E = \mathbf J^{-T} \tilde {\bm E}, \bm H = \mathbf J^{-T}\tilde {\bm H}$ and similarly the test functions $\tilde {\bm e},\tilde{\bm h}$ 
leads to the weak system
\begin{align*}
  \int_\Omega \left(\partial_t\det (\mathbf J)\mathbf J^{-1}\mathbf J^{-T} \tilde {\bm E}\right)\cdot \tilde {\bm e} &=\int_\Omega \curl \tilde {\bm H}\cdot \tilde {\bm e},\\
  \int_\Omega \left(\partial_t\det (\mathbf J)\mathbf J^{-1}\mathbf J^{-T} \tilde {\bm H}\right)\cdot \tilde {\bm h} &=-\int_\Omega \curl \tilde {\bm E}\cdot \tilde {\bm h}.
\end{align*}
For $\mathbf x\in\Omega_{pml}$ we have
\begin{align*}
  \partial_t \mathbf J^{-1}\mathbf J^{-T}\det (\mathbf J) &= \partial_t\frac{1}{1+\frac{\sigma}{\partial t}}\Pi_1+\partial_t(1+\frac{\sigma}{\partial t}))\Pi_1^\perp\\
  &=\partial_t \mathbf I  + \sigma (\Pi_1^\perp-\Pi_1)+\frac{\sigma^2}{\partial_t+\sigma}\Pi_1,
\end{align*}
Thus, introducing $\partial_t \hat{\bm E} = \sigma\tilde {\bm E}-\sigma\hat {\bm E}$ and
$\partial_t \hat {\bm H} = \sigma\tilde {\bm E}-\sigma\hat {\bm H}$, defined only on $\Omega_{pml}$
leads to the damped system
\begin{align}
  \begin{split}
  \int_\Omega \partial_t\tilde {\bm E}\cdot \tilde{\bm e} &=\sigma\int_{\Omega_{pml}} (\Pi_1-\Pi_1^\perp)\tilde {\bm E}\cdot \tilde{\bm e}-\sigma\int_{\Omega_{pml}}\Pi_1\hat {\bm E}\cdot \tilde{\bm e}+\int_\Omega \curl \tilde {\bm H}\cdot \tilde{\bm e},\\
  \int_{\Omega_{pml}} \partial_t\hat {\bm E}\cdot \hat{\bm e} &=\sigma\int_{\Omega_{pml}} \tilde {\bm E}\cdot \hat{\bm e}-\sigma\int_{\Omega_{pml}} \hat {\bm E}\cdot \hat{\bm e},\\
  \int_\Omega \partial_t\tilde {\bm H}\cdot \bm h &=\sigma\int_{\Omega_{pml}} (\Pi_1-\Pi_1^\perp)\tilde {\bm H}\cdot \tilde{\bm h}-\sigma\int_{\Omega_{pml}}\Pi_1\hat {\bm H}\cdot \tilde{\bm h}-\int_\Omega \curl \tilde {\bm E}\cdot \tilde{\bm h},\\
  \int_{\Omega_{pml}} \partial_t\hat {\bm H}\cdot \hat{\bm h} &=\sigma\int_{\Omega_{pml}} \tilde {\bm H}\cdot \hat{\bm h}-\sigma\int_{\Omega_{pml}} \hat {\bm H}\cdot \hat{\bm h}.\\
  \end{split}
  \label{eq:pml_system}
\end{align}
Now the damped system \eqref{eq:pml_system} is discretised in space using the discrete spaces from Section \ref{sec:spacedisc}, where the weak $\curl$ operators are discretised verbatim{, i.e., by replacing the continuous fields and test functions in \eqref{eq:pml_system} with their expansions in the primal and dual curl-conforming bases and by evaluating the resulting bilinear forms with the same quadrature rules as in the non-PML case.}

Similar to Section \ref{sec:timedisc}, we denote the fully discrete coefficient vectors of the discretisations of the fields $\tilde{\bm E}, \hat{\bm E},\tilde{\bm H},\hat{\bm H}$ at time $q\Delta t$ by $\tilde{\mathbf e}^q,\hat{\mathbf e}^q,\tilde{\mathbf h}^q,\hat {\mathbf h}^q$. Note that the discrete spaces for the auxiliary fields $\hat{\bm E}, \hat{\bm H}$ are defined merely on $\Omega_{pml}$ and thus have less components.
Using an extension of the leap frog scheme for damped equations we obtain the fully discrete system

\begin{align}
  \begin{split}
    & \tilde{\mathbf{h}}^{1/2} =\tilde{\mathbf{h}}^{0} +\frac{1}{2}\Delta t\mathbf{M}_{\boldsymbol{\mu}}^{-1} \left[\sigma\tilde{ \mathbf D}_{h}\tilde{\mathbf h}^0-\sigma\hat{\mathbf D}_{h}\hat{\mathbf h}^0- \mathbf{C} \tilde{\mathbf{e}}^0\right],\\
    &\hat{\mathbf h}^{1/2} = \hat{\mathbf h}^0+\frac{1}{2}\Delta t\sigma \left[\mathbf R_h\tilde{\mathbf h}^0-\hat{\mathbf h}^0\right],\\
    & \tilde{\mathbf{e}}^{q+1} =\tilde{\mathbf{e}}^{q} +\Delta t\mathbf{M}_{\boldsymbol{\varepsilon}}^{-1} \left[\sigma\tilde{\mathbf  D}_{e}\tilde{\mathbf e}^{q}-\sigma\hat{\mathbf D}_{e}\hat{\mathbf e}^{q}+ \mathbf{C}^T \tilde{\mathbf{h}}^{q+1/2}\right],\\
    &\hat{\mathbf e}^{q+1} = \hat{\mathbf e}^q+\Delta t\sigma \left[\mathbf R_e\tilde{\mathbf e}^q-\hat{\mathbf e}^q\right],\\
    & \tilde{\mathbf{h}}^{q+1/2} =\tilde{\mathbf{h}}^{q-1/2} +\Delta t\mathbf{M}_{\boldsymbol{\mu}}^{-1} \left[\sigma\tilde{\mathbf D}_{h}\tilde{\mathbf h}^{q-1/2}-\sigma\hat{\mathbf D}_{h}\hat{\mathbf h}^{q-1/2}- \mathbf{C} \tilde{\mathbf{e}}^q\right],\\
    &\hat{\mathbf h}^{q+1/2} = \hat{\mathbf h}^{q-1/2}+\Delta t\sigma \left[\mathbf R_h\tilde{\mathbf h}^{q-1/2}-\hat{\mathbf h}^{q-1/2}\right],\\
  \end{split}
  \label{eq:fully_discrete_pml}
\end{align}
where the entries of the damping matrices $\tilde{\mathbf D}_e,\hat{\mathbf D}_e,\tilde{\mathbf D}_h,\hat{\mathbf D}_h$ are defined similarly to \eqref{E9} by
\begin{align}
  \begin{split}
    &\tilde{D}_{e,mn}:=\sum_{K \in \mathcal{K},K\subset\Omega_{pml}} \int_{K}(\Pi_1-\Pi_1^\perp) \bm{E}_m\cdot \bm{E}_n\, d V,\\
    &\hat{D}_{e,mn}:=\sum_{K \in \mathcal{K},K\subset\Omega_{pml}} \int_{K}\Pi_1 \bm{E}_m\cdot \bm{E}_n\, d V,\\
    &\tilde{D}_{h,mn}:=\sum_{K \in \mathcal{K},K\subset\Omega_{pml}} \int_{K}(\Pi_1-\Pi_1^\perp) \bm{H}_m\cdot \bm{H}_n\, d V,\\
    &\hat{D}_{h,mn}:=\sum_{K \in \mathcal{K},K\subset\Omega_{pml}} \int_{K}\Pi_1 \bm{H}_m\cdot \bm{H}_n\, d V,
  \end{split}
\end{align}
and the matrices $\mathbf R_e,\mathbf R_h$ are the restriction operators from the whole discrete space defined on $\Omega$ to the subspace defined on $\Omega_{pml}$.

\subsection{Reference solution for the time-domain problem}
\label{sec:ref_sol}
We look for solutions of the Maxwell system on the domain
\begin{align*}
  \Omega = (0,l_x)\times(0,l_y)\times(0,l_z),
\end{align*}
with homogeneous boundary conditions on the sidewall, the incoming wave condition \eqref{eq:incoming_wave} at $\Gamma_{in}=\{0\}\times(0,l_y)\times(0,l_z)$, and no reflections at
$\Gamma_{out}=\{l_x\}\times(0,l_y)\times(0,l_z)$

We make an ansatz (
note that $\bm E$ fulfills the homogeneous boundary condition on the side walls)
\begin{align*}
\bm E(t,x,y,z)&=\left(0,0,E_z\right)^\top=\left(0,0,e(t,x)\sin(y\pi l_y^{-1})\right)^\top,\\
\bm H(t,x,y,z)&=\left(H_x,H_y,0\right)^\top=\left(h_x(t,x)\cos(y\pi l_y^{-1}),h_y(t,x)\sin(y\pi l_y^{-1}),0\right)^\top.
\end{align*}

This immediately gives
\begin{align*}
\mathrm{curl}\mathrm{curl}\bm E=\mathrm{curl}\left(\partial_yE_z,0,-\partial_x E_z\right)^\top=\left(0,0,-(\partial^2_x+\partial_y^2)E_z\right)^\top=\left(0,0,(\pi^2l_y^{-2}-\partial_x^2)e(t,x)\sin(y\pi l_y^{-1})\right)^\top.
\end{align*}
Thus $e$ solves the one dimensional problem
\begin{align*}
\partial_t^2e(t,x)&=(\partial_x^2-\pi^2l_y^{-2})e(t,x),\\
e(0,x)&=0,\\
e(0,l_x)&\text{ is propagating},\\
e(t,0)&=e_0(t).
\end{align*}
Above equation can be solved analytically in Laplace domain and using known Laplace transforms it can be written as a convolution in the following way:
\begin{align*}
e(t,x)=\begin{cases}
e_0(t-x)-\left(e_0(\cdot)*\frac{\pi}{l_y}x\frac{J_1\left(\frac{\pi}{l_y}\sqrt{\cdot^2-x^2}\right)}{\sqrt{\cdot^2-x^2}}\Theta(\cdot-x)\right)(t),&t\geq x,\\
0,&t<x.
\end{cases}
\end{align*}
{Here $J_1$ denotes the Bessel function of the first kind of order one.}
Assuming that $e_0$ vanishes for negative times the convolution for $t>x$ is given by the integral, 
\begin{align*}
  \left(e_0(\cdot)*\frac{\pi}{l_y}x\frac{J_1\left(\frac{\pi}{l_y}\sqrt{\cdot^2-x^2}\right)}{\sqrt{\cdot^2-x^2}}\Theta(\cdot-x)\right)(t)
  {}&=\frac{\pi}{l_y}x\int_{-\infty}^\infty e_0(t-\tau)\frac{J_1\left(\frac{\pi}{l_y}\sqrt{\tau^2-x^2}\right)}{\sqrt{\tau^2-x^2}}\Theta(\tau-x) d\tau\\
  {}&=\frac{\pi}{l_y}x\int_{x}^{t} e_0(t-\tau)\frac{J_1\left(\frac{\pi}{l_y}\sqrt{\tau^2-x^2}\right)}{\sqrt{\tau^2-x^2}} d\tau,\\
\end{align*}
and can be approximated using a suitable numerical quadrature rule.

\bibliographystyle{unsrt}
\bibliography{bibliography}

 \end{document}